\documentclass[a4paper,12pt]{article}
\usepackage[pagewise]{lineno}

\newcommand*\patchAmsMathEnvironmentForLineno[1]{
  \expandafter\let\csname old#1\expandafter\endcsname\csname #1\endcsname
  \expandafter\let\csname oldend#1\expandafter\endcsname\csname end#1\endcsname
  \renewenvironment{#1}
     {\linenomath\csname old#1\endcsname}
     {\csname oldend#1\endcsname\endlinenomath}}
\newcommand*\patchBothAmsMathEnvironmentsForLineno[1]{
  \patchAmsMathEnvironmentForLineno{#1}
  \patchAmsMathEnvironmentForLineno{#1*}}
\AtBeginDocument{
\patchBothAmsMathEnvironmentsForLineno{equation}
\patchBothAmsMathEnvironmentsForLineno{align}
\patchBothAmsMathEnvironmentsForLineno{flalign}
\patchBothAmsMathEnvironmentsForLineno{alignat}
\patchBothAmsMathEnvironmentsForLineno{gather}
\patchBothAmsMathEnvironmentsForLineno{multline}
}

\usepackage{amsmath}
\usepackage{latexsym}
\usepackage{url}
\usepackage{color}
\numberwithin{equation}{section}

\def\R{{\bf R}}

\def\N{{\bf N}}

\def\d{\displaystyle}
\def\e{{\varepsilon}}

\def\wt{\widetilde}

\def\p{\partial}
\def\v#1{\mbox{\boldmath $#1$}}

\newtheorem{thm}{Theorem}[section]

\newtheorem{prop}{Proposition}[section]
\newtheorem{rem}{Remark}[section]

\title{The combined effect in one space dimension beyond the general theory\\
for nonlinear wave equations\\
{\small\it Dedicated to Professor Tohru Ozawa on his sixtieth birthday}}
\author{
Katsuaki Morisawa\\
{\footnotesize Musashi High School and Junior High School,}\\
{\footnotesize1-26-1 Toyotamakami, Nerima, Tokyo, 176-8535, Japan.,}\\
\\
Takiko Sasaki\\
{\footnotesize
Department of Mathematical Engineering, Faculty of Engineering, Musashino University,}\\
{\footnotesize
3-3-3 Ariake, Koto-ku, Tokyo 135-8181, Japan.}\\
{\footnotesize /Mathematical Institute, Tohoku University,
Aoba, Sendai 980-8578, Japan.,}
\\
\\
Hiroyuki Takamura
\\
{\footnotesize Mathematical Institute, Tohoku University,}\\
{\footnotesize Aoba, Sendai 980-8578, Japan.}
\date{}
}
\begin{document}
\maketitle
\begin{abstract}
In this paper, we show the so-called \lq\lq combined effect" of two different kinds
of nonlinear terms for semilinear wave equations in one space dimension.
Such a special phenomenon appears only in case
the total integral of the initial speed is zero.
It is remarkable that,
including the combined effect case, our results on the lifespan estimates
are partially better than those of the general theory for nonlinear wave equations.  
\end{abstract}


\section{Introduction}

\par 
We consider the initial value problems;
\begin{equation}
\label{IVP_combined}
\left\{
\begin{array}{ll}
	\d u_{tt}-u_{xx}=A|u_t|^p+B|u|^q
	&\mbox{in}\quad \R\times(0,\infty),\\
	u(x,0)=\e f(x),\ u_t(x,0)=\e g(x),
	& x\in\R,
\end{array}
\right.
\end{equation}
where $p,q>1$, $A,B\ge0$, $f$ and $g$ are given smooth functions of compact support
and a parameter $\e>0$ is \lq\lq small enough".
We are interested in the lifespan $T(\e)$, the maximal existence time,
of classical solutions of (\ref{IVP_combined}).
Our results in this paper are the following estimates for $A>0$ and $B>0$:
\begin{equation}
\label{lifespan_non-zero}
T(\e)\sim
\min\{C\e^{-(p-1)},C\e^{-(q-1)/2}\}\qquad \mbox{if}\ \int_{\R}g(x)dx\not=0\\
\end{equation}
and
\begin{equation}
\label{lifespan_zero}
\begin{array}{c}
T(\e)\sim
\left\{
\begin{array}{ll}
C\e^{-p(q-1)/(q+1)} & \mbox{for}\ \d\frac{q+1}{2}\le p\le q,\\
\min\{C\e^{-(p-1)},C\e^{-q(q-1)/(q+1)}\} & \mbox{otherwise}\\
\end{array}
\right.
\\
 \mbox{if}\ \d\int_{\R}g(x)dx=0.
 \end{array}
\end{equation}
Here we denote the fact that there are positive constants,
$C_1$ and $C_2$, independent of $\e$ satisfying $A(\e,C_1)\le T(\e)\le A(\e,C_2)$
by $T(\e)\sim A(\e,C)$.

\par
Recall that we have
\[
T(\e)\sim C\e^{-(p-1)}\quad\mbox{for $A>0$ and $B=0$.}
\]
This result was verified by Zhou \cite{Zhou01} for the upper bound,
and by Li,Yu and Zhou \cite{LYZ91,LYZ92}
for the lower bound with integer $p\ge2$ including more general nonlinear term.
We call \cite{LYZ91,LYZ92} \lq\lq general theory" for nonlinear wave equations in one dimension.
Recently, Kitamura, Morisawa and Takamura \cite{KMT} have verified the lower bound
for all $p\ge1$ including the case that nonlinear term has spatial weights,
in which only the $C^1$ solution of the associated integral equation is considered for $1<p<2$.
But it can be also the classical solution by trivial modifications on estimating 
the nonlinear term with H\"older continuity.
See Remark \ref{rem:regularity} below. 
On the other hand, Zhou \cite{Zhou92} obtained
\[
T(\e)\sim
\left\{
\begin{array}{ll}
C\e^{-(q-1)/2} & \mbox{if}\ \d\int_{\R}g(x)dx\neq0,\\
C\e^{-q(q-1)/(q+1)} & \mbox{if}\ \d\int_{\R}g(x)dx=0
\end{array}
\right.
\mbox{for $A=0$ and $B>0$}.
\]
Therefore (\ref{lifespan_non-zero}) and (\ref{lifespan_zero}) are quite natural
as taking the minimum of both results
except for the first case in (\ref{lifespan_zero}),
in which, we have
\[
C\e^{-p(q-1)/(q+1)}\le\min\{C\e^{-(p-1)},C\e^{-q(q-1)/(q+1)}\}
\quad\mbox{for}\ \d\frac{q+1}{2}\le p\le q.
\]
We call this special phenomenon by \lq\lq combined effect" of two nonlinearities.
The combined effect was first observed by
Han and Zhou in \cite{HZ14}, which targets to show the optimality of the result
of Katayama \cite{Katayama01} on the lower bound of the lifespan of classical solutions
of nonlinear wave equations with a nonlinear term $u_t^3+u^4$ in two space dimensions
including more general nonlinear terms.
It is known that $T(\e)\sim\exp\left(C\e^{-2}\right)$ for the nonlinear term
$u_t^3$ and $T(\e)=\infty$ for the nonlinear term $u^4$,
but Katayama \cite{Katayama01} obtained only a much worse estimate
than their minimum  as $T(\e)\ge c\e^{-18}$.
Surprisingly, more than ten years later, Han and Zhou \cite{HZ14} showed that
this result is optimal as $T(\e)\le C\e^{-18}$.
They also considered (\ref{IVP_combined}) for all space dimensions $n$ bigger than 1
and obtain the upper bound of the lifespan.
Its counterpart, the lower bound of the lifespan, was obtained by
Hidano, Wang and Yokoyama \cite{HWY16} for $n=2,3$.
See the introduction of \cite{HWY16} for the precise results and references.  
We note that the first case in (\ref{lifespan_zero}) coincides with the lifespan estimate 
for the combined effect in \cite{HZ14, HWY16} if one sets $n=1$ formally.
Later, Dai, Fang and Wang \cite{DFW19} improved the lower bound of lifespan
for the critical case in Hidano, Wang and Yokoyama \cite{HWY16}.
They also showed that $T(\e)<\infty$ for all $p,q>1$ in case of $n=1$, i.e. (\ref{IVP_combined}).

\par
Finally we strongly remark that our estimates in (\ref{lifespan_non-zero}) and (\ref{lifespan_zero})
are better than those of the general theory by Li, Yu and Zhou \cite{LYZ91, LYZ92}:
\[
T(\e)\ge
\left\{
\begin{array}{ll}
C\e^{-(p-1)/2} & \mbox{if}\ \d\int_{\R}g(x)dx\neq0,\\
C\e^{-p(p-1)/(p+1)} & \mbox{if}\ \d\int_{\R}g(x)dx=0
\end{array}
\right.\]
in case of
\[
\frac{q+1}{2}<p<q
\]
with integer $p,q\ge2$.
Indeed, our result on the lower bound of the lifespan can be established
also for the indefinite sign  terms as $u_{tt}-u_{xx}=u_t^p+u^q$. 
The typical example is $(p,q)=(4,5)$ in such a case.
This fact shows a possibility to improve the general theory.
For details, see the last half of the next section.
Of course, some special structure of the nonlinear terms
such as \lq\lq null condition" guarantees the global-in-time existence.
See Nakamura \cite{Nakamura14}, Luli, Yang and Yu \cite{LYY18},
Zha \cite{Zha20, Zha22} for examples in this direction.
But we are interested in the optimality of the general theory.
The details are discussed at the end of Section 2 below.

\par
This paper is organized as follows.
In the next section, the preliminaries are introduced.
Moreover, (\ref{lifespan_non-zero}) and (\ref{lifespan_zero}) are divided into four theorems,
and we compare our results with those of the general theory.
Sections 3 and 4 are devoted to the proof of the existence part of (\ref{lifespan_non-zero}).
Sections 5, 6 and 7 are devoted to the proof of the existence part of (\ref{lifespan_zero}).
Their main strategy is the iteration method in the weighted $L^\infty$ space 
due to Kitamura, Morisawa and Takamura \cite{KMT22, KMT},
which is originally introduced by John \cite{John79}.
Finally, we prove the blow-up part of (\ref{lifespan_non-zero}) and (\ref{lifespan_zero}),
which is different from the functional method by Han and Zhou \cite{HZ14}.


\section{Preliminaries and main results}

Throughout this paper, we assume that the initial data
$(f,g)\in C_0^2(\R)\times C^1_0(\R)$ satisfies
\begin{equation}
\label{supp_initial}
\mbox{\rm supp }f,\ \mbox{supp }g\subset\{x\in\R:|x|\le R\},\quad R\ge1.
\end{equation}
Let $u$ be a classical solution of (\ref{IVP_combined}) in the time interval $[0,T]$.
Then the support condition of the initial data, (\ref{supp_initial}), implies that
\begin{equation}
\label{support_sol}
\mbox{supp}\ u(x,t)\subset\{(x,t)\in\R\times[0,T]:|x|\le t+R\}.
\end{equation}
For example, see Appendix of John \cite{John_book} for this fact.

\par
It is well-known that $u$ satisfies the following integral equation.
\begin{equation}
\label{u}
u(x,t)=\e u^0(x,t)+L(A|u_t|^p+B|u|^q)(x,t),
\end{equation}
where $u^0$ is a solution of the free wave equation with the same initial data,
\begin{equation}
\label{u^0}
u^0(x,t):=\frac{1}{2}\{f(x+t)+f(x-t)\}+\frac{1}{2}\int_{x-t}^{x+t}g(y)dy,
\end{equation}
and a linear integral operator $L$ for a function $v=v(x,t)$ in Duhamel's term is defined by
\begin{equation}
\label{nonlinear}
L(v)(x,t):=\frac{1}{2}\int_0^tds\int_{x-t+s}^{x+t-s}v(y,s)dy.
\end{equation}
Then, one can apply the time-derivative to (\ref{u}) to obtain
\begin{equation}
\label{u_t}
u_t(x,t)=\e u_t^0(x,t)+L'(A|u_t|^p+B|u|^q)(x,t)
\end{equation}
and
\begin{equation}
\label{u^0_t}
u_t^0(x,t)=\frac{1}{2}\{f'(x+t)-f'(x-t)+g(x+t)+g(x-t)\},
\end{equation}
where $L'$ for a function $v=v(x,t)$ is defined by
\begin{equation}
\label{nonlinear_derivative}
L'(v)(x,t):=\frac{1}{2}\int_0^t\{v(x+t-s,s)+v(x-t+s,s)\}ds.
\end{equation}
On the other hand, applying the space-derivative to (\ref{u}),
we have
\begin{equation}
\label{u_x}
u_x(x,t)=\e u_x^0(x,t)+\overline{L'}(A|u_t|^p+B|u|^q)(x,t)
\end{equation}
and
\begin{equation}
\label{u^0_x}
u_x^0(x,t)=\frac{1}{2}\{f'(x+t)+f'(x-t)+g(x+t)-g(x-t)\},
\end{equation}
where $\overline{L'}$ for a function $v=v(x,t)$ is defined by
\begin{equation}
\label{nonlinear_derivative_conjugate}
\overline{L'}(v)(x,t):=
\frac{1}{2}\int_0^t\{v(x+t-s,s)-v(x-t+s,s)\}ds.
\end{equation}
Therefore, $u_x$ is expressed by $u$ and $u_t$.
Moreover, applying one more space-derivative to (\ref{u_t}) yields that
\begin{equation}
\label{u_tx}
\begin{array}{ll}
u_{tx}(x,t)=&\d\e u_{tx}^0(x,t)\\
&\d+L'(Ap|u_t|^{p-2}u_tu_{tx}+Bq|u|^{q-2}uu_x)(x,t)
\end{array}
\end{equation}
and
\begin{equation}
\label{u^0_tx}
u_{tx}^0(x,t)=\frac{1}{2}\{f''(x+t)-f''(x-t)+g'(x+t)+g'(x-t)\}.
\end{equation}
Similarly, we have that
\[
\begin{array}{ll}
u_{tt}(x,t)=&
\d\e u_{tt}^0(x,t)+A|u_t(x,t)|^p+B|u(x,t)|^q\\
&\d +\overline{L'}(Ap|u_t|^{p-2}u_tu_{tx}+Bq|u|^{q-2}uu_x)(x,t)
\end{array}
\]
and
\[
u_{tt}^0(x,t)=\frac{1}{2}\{f''(x+t)+f''(x-t)+g'(x+t)-g'(x-t)\}.
\]
Therefore, $u_{tt}$ is expressed by $u,u_t,u_x,u_{tx}$ and so is $u_{xx}$ because of
\[
\begin{array}{ll}
u_{xx}(x,t)=
&\d\e u_{xx}^0(x,t)\\
&\d+\overline{L'}(Ap|u_t|^{p-2}u_tu_{tx}+Bq|u|^{q-2}uu_x)(x,t)
\end{array}
\]
and
\[
u_{xx}^0(x,t)=u^0_{tt}(x,t).
\]

\begin{rem}
\label{rem:regularity}
In view of (\ref{u_tx}),
it is sufficient to employ H\"older continuity of the nonlinear term, i.e.
\[
\left||a|^{p-2}a-|b|^{p-2}b\right|\le2|a-b|^{p-1}\ (a,b\in\R,\ 1<p<2),
\]
in estimating the difference of the approximation sequence
to construct a classical solution for $p,q\in(1,2)$.
This fact was overlooked in Kitamura, Morisawa and Takamura \cite{KMT}
as stated in Introduction.
\end{rem}

\par
First, we note the following fact.

\begin{prop}
\label{prop:system}
Assume that $(f,g)\in C^2(\R)\times C^1(\R)$.
Let $(u,w)$ be a $C^1$ solution of a system of integral equations
\begin{equation}
\label{system}
\left\{
\begin{array}{l}
u=\e u^0+L(A|w|^p+B|u|^q),\\
w=\e u_t^0+L'(A|w|^p+B|u|^q)
\end{array}
\right.
\mbox{in}\ \R\times[0,T]
\end{equation}
with some $T>0$.
Then, $w\equiv u_t$ in $\R\times[0,T]$ holds
and $u$ is a classical solution of (\ref{IVP_combined}) in $\R\times[0,T]$.
\end{prop}
\par\noindent
{\bf Proof.} It is trivial that $w\equiv u_t$ by differentiating the first equation with respect to $t$.
 The rest part is easy along with the computation above in this section. 
\hfill$\Box$

\vskip10pt
Our results in (\ref{lifespan_non-zero}) and (\ref{lifespan_zero})
are divided into the following four theorems.

\begin{thm}
\label{thm:lower-bound_non-zero}
Let $A>0$ and $B>0$.
Assume (\ref{supp_initial}) and
\begin{equation}
\label{non-zero}
\int_{\R}g(x)dx\neq0.
\end{equation}
Then, there exists a positive constant $\e_1=\e_1(f,g,p,q,A,B,R)>0$ such that
a classical solution $u\in C^2(\R\times[0,T])$ of (\ref{IVP_combined}) exists
as far as $T$ satisfies
\begin{equation}
\label{lower-bound_non-zero}
T\le
\left\{
\begin{array}{ll}
c\e^{-(p-1)} & \mbox{for}\ p\le\d\frac{q+1}{2},\\
c\e^{-(q-1)/2} & \mbox{for}\ \d\frac{q+1}{2}\le p,
\end{array}
\right.
\end{equation}
where $0<\e\le\e_1$ and $c$ is a positive constant independent of $\e$.
\end{thm}

\begin{thm}
\label{thm:lower-bound_zero}
Let $A>0$ and $B>0$.
Assume (\ref{supp_initial}) and
\begin{equation}
\label{zero}
\int_{\R}g(x)dx=0.
\end{equation}
Then, there exists a positive constant $\e_2=\e_2(f,g,p,q,A,B,R)>0$ such that
a classical solution $u\in C^2(\R\times[0,T])$ of (\ref{IVP_combined}) exists
as far as $T$ satisfies
\begin{equation}
\label{lower-bound_zero}
T\le
\left\{
\begin{array}{ll}
c\e^{-(p-1)} & \mbox{for}\ p\le\dfrac{q+1}{2},\\
c\e^{-p(q-1)/(q+1)} & \mbox{for}\ \d\frac{q+1}{2}\le p\le q,\\
c\e^{-q(q-1)/(q+1)} & \mbox{for}\ p\ge q,
\end{array}
\right.
\end{equation}
where $0<\e\le\e_2$ and $c$ is a positive constant independent of $\e$.
\end{thm}

\begin{thm}
\label{thm:upper-bound_non-zero}
Let $A>0$ and $B>0$.
Assume (\ref{supp_initial}) and
\begin{equation}
\label{positive_non-zero}
\int_{\R}g(x)dx>0.
\end{equation}
Then, there exists a positive constant $\e_3=\e_3(f,g,p,q,A,B,R)>0$ such that
a classical solution of (\ref{IVP_combined}) cannot exist
as far as $T$ satisfies
\begin{equation}
\label{upper-bound_non-zero}
T\ge
\left\{
\begin{array}{ll}
C\e^{-(p-1)} & \mbox{for}\ p\le\d\frac{q+1}{2},\\
C\e^{-(q-1)/2} & \mbox{for}\ \d\frac{q+1}{2}\le p,
\end{array}
\right.
\end{equation}
where $0<\e\le\e_3$ and $C$ is a positive constant independent of $\e$.
\end{thm}

\begin{thm}
\label{thm:upper-bound_zero}
Let $A>0$ and $B>0$.
Assume (\ref{supp_initial}) and
\begin{equation}
\label{positive_zero}
f(x)\ge0(\not\equiv0),\ f'(x)<0\ \mbox{for}\ x\in(0,R)\ \mbox{and}\ g(x)\equiv0.
\end{equation}
Then, there exists a positive constant $\e_4=\e_4(f,p,q,A,B,R)>0$ such that
a classical solution of (\ref{IVP_combined}) cannot exist
as far as $T$ satisfies
\begin{equation}
\label{upper-bound_zero}
T\ge
\left\{
\begin{array}{ll}
C\e^{-(p-1)} & \mbox{for}\ p\le\dfrac{q+1}{2},\\
C\e^{-p(q-1)/(q+1)} & \mbox{for}\ \d\frac{q+1}{2}\le p\le q,\\
C\e^{-q(q-1)/(q+1)} & \mbox{for}\ p\ge q,
\end{array}
\right.
\end{equation}
where $0<\e\le\e_4$ and $C$ is a positive constant independent of $\e$.
\end{thm}

\begin{rem}
\label{rem:relation}
It is trivial that Theorem \ref{thm:lower-bound_non-zero} and
Theorem \ref{thm:upper-bound_non-zero} imply (\ref{lifespan_non-zero}).
On the other hand, we have that
\[
p-1=\frac{q(q-1)}{q+1}\quad\Longleftrightarrow\quad p=\frac{q^2+1}{q+1}
\]
and
\[
\frac{q+1}{2}<\frac{q^2+1}{q+1}<q.
\]
Moreover, we see that
\[
p-1\le\frac{p(q-1)}{q+1}\quad\Longleftrightarrow\quad p\le\frac{q+1}{2}.
\]
Therefore Theorem \ref{thm:lower-bound_zero} and
Theorem \ref{thm:upper-bound_zero} imply (\ref{lifespan_zero}).
\end{rem}

\par
The proofs of four theorems above appear in the following sections.
Form now on, we shall compare our results with those of the general theory 
by Li, Yu and Zhou \cite{LYZ91, LYZ92},
in which the following problem of general form is considered:
\begin{equation}
\label{IVP_general}
\left\{
\begin{array}{ll}
	\d u_{tt}-u_{xx}=F(u,Du,\p_xDu)
	&\mbox{in}\quad \R\times(0,\infty),\\
	u(x,0)=\e f(x),\ u_t(x,0)=\e g(x),
	& x\in\R,
\end{array}
\right.
\end{equation}
where we denote $D:=(\p_t,\p_x)$ and $F\in C^\infty(\R^5)$ satisfies
\[
F(\lambda)=O(|\lambda|^{1+\alpha})\quad\mbox{with $\alpha\in\N$ near $\lambda=0$}.
\]
(\ref{IVP_general}) requires $f,g\in C_0^\infty(\R)$.
Then, the lifespan of the classical solution of (\ref{IVP_general}) defined by $\wt{T}(\e)$ has
estimates from below as
\begin{equation}
\label{lifespan_general}
\wt{T}(\e)\ge
\left\{
\begin{array}{ll}
c\e^{-\alpha/2} & \mbox{in general},\\
c\e^{-\alpha(1+\alpha)/(2+\alpha)} & \mbox{if}\ \d\int_{\R}g(x)dx=0,\\
c\e^{-\alpha} & \mbox{if $\p_u^\beta F(0)=0$ for $1+\alpha\le\forall\beta\le2\alpha$}.
\end{array}
\right.
\end{equation}
This is the result of the general theory.
If one applies it to our problem (\ref{IVP_combined}) with
\begin{equation}
\label{F_special}
F(u,Du,\p_xDu)=u_t^p+u^q\quad\mbox{with}\ 2\le p,q\in\N,
\end{equation}
one has the following estimates in each cases.
\begin{itemize}
\item
When $p<q$,
\par
then, we have to set $\alpha=p-1$ which yields that
\[
\wt{T}(\e)\ge
\left\{
\begin{array}{ll}
c\e^{-(p-1)/2} & \mbox{in general},\\
c\e^{-p(p-1)/(p+1)} & \mbox{if}\ \d\int_{\R}g(x)dx=0,\\
c\e^{-(p-1)} & \mbox{if $\p_u^\beta F(0)=0$ for $p\le\forall\beta\le2(p-1)$}.
\end{array}
\right.
\]
For $p\le(q+1)/2$, i.e. $2p-1\le q$, we see that
\[
\p_u^\beta F(u,Du,\p_xDu)=O(u^{q-\beta})
\quad\mbox{and}\quad 1\le2p-1-\beta\le q-\beta
\]
because of $p+1\le2p-1$, which yields
\[
\p_u^\beta F(0)=0\quad\mbox{for}\ p\le\forall\beta\le2(p-1).
\]
Therefore the third case holds and we obtain
\[
\wt{T}(\e)\ge c\e^{-(p-1)}
\] 
whatever the value of $\d\int_{\R}g(x)dx$ is.
On the other hand, for $(q+1)/2<p$, i.e. $(p<)q<2p-1$, we see that
\[
\exists\beta\in\{p+1,\ldots,2p-2\}\mbox{ s.t. }\p_u^\beta F(0)\neq0,
\]
so that the third case does not hold and we obtain
\[
\wt{T}(\e)\ge
\left\{
\begin{array}{ll}
c\e^{-(p-1)/2} & \mbox{if}\ \d\int_{\R}g(x)dx\neq0,\\
c\e^{-p(p-1)/(p+1)} & \mbox{if}\ \d\int_{\R}g(x)dx=0.\\
\end{array}
\right.
\]
\item
When $p\ge q$,
\par
then, similarly to the case above, we have to set $\alpha=q-1$, which yields that
\[
\wt{T}(\e)\ge
\left\{
\begin{array}{ll}
c\e^{-(q-1)/2} & \mbox{in general},\\
c\e^{-q(q-1)/(q+1)} & \mbox{if}\ \d\int_{\R}g(x)dx=0,\\
c\e^{-(q-1)} & \mbox{if $\p_u^\beta F(0)=0$ for $q\le\forall\beta\le2(q-1)$}.
\end{array}
\right.
\]
We note that the third case does not hold by $\p_u^q F(0)\neq0$.
\end{itemize}
In conclusion, for the special nonlinear term in (\ref{F_special}),
the result of the general theory is
\[
\wt{T}(\e)\ge
\left\{
\begin{array}{ll}
c\e^{-(p-1)} & \mbox{for}\ p\le\d\frac{q+1}{2},\\
c\e^{-(p-1)/2} & \mbox{for}\ \d\frac{q+1}{2}\le p\le q,\\
c\e^{-(q-1)/2} & \mbox{for}\ q\le p
\end{array}
\right.
\quad\mbox{if}\ \int_{\R}g(x)dx\not=0\\
\]
and
\[
\wt{T}(\e)\ge
\left\{
\begin{array}{ll}
c\e^{-(p-1)} & \mbox{for}\ p\le\d\frac{q+1}{2},\\
c\e^{-p(p-1)/(p+1)} & \mbox{for}\ \d\frac{q+1}{2}\le p\le q,\\
c\e^{-q(q-1)/(q+1)} & \mbox{for}\ q\le p\\
\end{array}
\right.
\quad \mbox{if}\ \d\int_{\R}g(x)dx=0.
\]
Therefore a part of our results,
\begin{equation}
\label{better}
T(\e)\sim
\left\{
\begin{array}{ll}
C\e^{-(q-1)/2} & \mbox{if}\ \d\int_{\R}g(x)\neq0,\\
C\e^{-p(q-1)/(q+1)} & \mbox{if}\ \d\int_{\R}g(x)dx=0
\end{array}
\right.
\quad\mbox{for}\ \frac{q+1}{2}\le p\le q,
\end{equation}
is better than the lower bound of $\wt{T}(\e)$.
If one follows the proof in the following sections, one can find that it is easy to see that
our results on the lower bounds also hold for a special term (\ref{F_special})
by estimating the difference of nonlinear terms from above
after employing the mean value theorem.
We note that we have infinitely many examples of $(p,q)=(m,m+1)$ as the inequality
\[
\frac{q+1}{2}=\frac{m+2}{2}<p=m<q=m+1
\]
holds for $m=3,4,5,\ldots$.
This fact indicates that we still have a possibility to improve the general theory
in the sense that the optimal results in (\ref{better}) should be included at least.


\section{Proof of Theorem \ref{thm:lower-bound_non-zero}}
\par
According to Proposition \ref{prop:system},
we shall construct a $C^1$ solution of (\ref{system}).
Let $\{(u_j,w_j)\}_{j\in\N}$ be a sequence of $\{C^1(\R\times[0,T])\}^2$ defined by
\begin{equation}
\label{u_j,w_j}
\left\{
\begin{array}{ll}
u_{j+1}=\e u^0+L(A|w_j|^p+B|u_j|^q), & u_1=\e u^0,\\
w_{j+1}=\e u_t^0+L'(A|w_j|^p+B|u_j|^q), &  w_1=\e u_t^0.
\end{array}
\right.
\end{equation}
Then, in view of (\ref{u_x}) and (\ref{u_tx}), $\left((u_j)_x,(w_j)_x\right)$ has to satisfy
\begin{equation}
\label{u_j_x,w_j_x}
\left\{
\begin{array}{l}
(u_{j+1})_x=\e u^0_x+\overline{L'}\left(A|w_j|^p+B|u_j|^q\right),\\
(u_1)_x=\e u^0_x,\\
(w_{j+1})_x=\e u^0_{tx}+L'\left(Ap|w_j|^{p-2}w_j(w_j)_x+Bq|u_j|^{q-2}u_j(u_j)_x\right),\\
 (w_1)_x=\e u^0_{tx},
\end{array}
\right.
\end{equation}
so that the function space in which $\{(u_j,w_j)\}$ converges is
\[
\begin{array}{ll}
X:=&\{(u,w)\in\{C^1(\R\times[0,T])\}^2\ :\ \|(u,w)\|_X<\infty,\\
&\quad\mbox{supp}\ (u,w)\subset\{(x,t)\in\R\times[0,T]\ :\ |x|\le t+R\}\},
\end{array}
\]
which is equipped with a norm
\[
\|(u,w)\|_X:=\|u\|_1+\|u_x\|_1+\|w\|_2+\|w_x\|_2,
\]
where
\[
\begin{array}{l}
\d\|u\|_1:=\sup_{(x,t)\in\R\times[0,T]}|u(x,t)|,\\
\d\|w\|_2:=\sup_{(x,t)\in\R\times[0,T]}|(t-|x|+2R)w(x,t)|.
\end{array}
\]
First we note that supp $(u_j,w_j)\subset\{(x,t)\in\R\times[0,T]\ :\ |x|\le t+R\}$ implies supp
$(u_{j+1},w_{j+1})\subset\{(x,t)\in\R\times[0,T]\ :\ |x|\le t+R\}$.
It is easy to check this fact by assumption on the initial data (\ref{supp_initial})
and the definitions of $L,\overline{L},L',\overline{L'}$ in the previous section.

\par
The following lemma contains some useful a priori estimates.
\begin{prop}
\label{prop:apriori}
Let $(u,w)\in\{C(\R\times[0,T])\}^2$ and supp\ $(u,w)\subset\{(x,t)\in\R\times[0,T]:|x|\le t+R\}$. Then there exists a positive constant $C$ independent of $T$ and $\e$ such that
\begin{equation}
\label{apriori}
\begin{array}{l}
\|L(|w|^p)\|_1\le C\|w\|_2^p(T+R),\ \|L(|u|^q)\|_1\le C\|u\|_1^q(T+R)^2,\\
\|L'(|w|^p)\|_2\le C\|w\|_2^p(T+R),\ \|L'(|u|^q)\|_2\le C\|u\|_1^q(T+R)^2,\\
\|L'(|w|^p)\|_1\le C\|w\|_2^p(T+R),\ \|L'(|u|^q)\|_1\le C\|u\|_1^q(T+R)^2.
\end{array}
\end{equation}
\end{prop}
The proof of Proposition \ref{prop:apriori} is established in the next section.
Set
\[
M:=\sum_{\alpha=0}^2\|f^{(\alpha)}\|_{L^\infty(\R)}
+\|g\|_{L^1(\R)}+\sum_{\beta=0}^1\|g^{(\beta)}\|_{L^\infty(\R)}.
\]
\vskip10pt
\par\noindent
{\bf The convergence of the sequence $\v{\{(u_j,w_j)\}}$.}
\par
First we note that $\|u_1\|_1,\|w_1\|_2\le M\e$ by (\ref{u^0}) and (\ref{u^0_t}).
Since (\ref{u_j,w_j}) and (\ref{apriori}) yield that
\[
\left\{
\begin{array}{ll}
\|u_{j+1}\|_1
&\le M\e+A\|L(|w_j|^p)\|_1+B\|L(|u_j|^q))\|_1\\
&\le M\e+AC\|w_j\|_2^p(T+R)+BC\|u_j\|^q_1(T+R)^2,\\
\|w_{j+1}\|_2
&\le M\e+A\|L'(|w_j|^p)\|_2+B\|L'(|u_j|^q)\|_2\\
&\le M\e+AC\|w_j\|_2^p(T+R)+CB\|u_j\|^q_1(T+R)^2,\
\end{array}
\right.
\]
the boundedness of $\{(u_j,w_j)\}$, i.e.
\begin{equation}
\label{bound_(u,w)}
\|u_j\|_1,\|w_j\|_2\le 3M\e\quad(j\in\N),
\end{equation}
follows from
\begin{equation}
\label{condi1}
AC(3M\e)^p(T+R),BC(3M\e)^q(T+R)^2\le M\e.
\end{equation}
Assuming (\ref{condi1}), one can estimate $(u_{j+1}-u_j)$ and $(w_{j+1}-w_j)$ as follows.
\[
\begin{array}{ll}
\|u_{j+1}-u_j\|_1
&\le\|L(A|w_j|^p-A|w_{j-1}|^p+B|u_j|^q-B|u_{j-1}|^q)\|_1\\
&\le 2^{p-1}pA\|L\left((|w_j|^{p-1}+|w_{j-1}|^{p-1})|w_j-w_{j-1}|\right)\|_1\\
&\quad+2^{q-1}qB\|L\left((|u_j|^{q-1}+|u_{j-1}|^{q-1})|u_j-u_{j-1}|\right)\|_1\\
&\le 2^{p-1}pAC(T+R)(\|w_j\|_2^{p-1}+\|w_{j-1}\|_2^{p-1})\|w_j-w_{j-1}\|_2\\
&\quad+2^{q-1}qBC(T+R)^2(\|u_j\|_1^{q-1}+\|u_{j-1}\|_1^{q-1})\|u_j-u_{j-1}\|_1\\
&\le 2^ppAC(3M\e)^{p-1}(T+R)\|w_j-w_{j-1}\|_2\\
&\quad+2^qqBC(3M\e)^{q-1}(T+R)^2\|u_j-u_{j-1}\|_1
\end{array}
\]
and
\[
\begin{array}{ll}
\|w_{j+1}-w_j\|_2
&\le\|L'(A|w_j|^p-A|w_{j-1}|^p+B|u_j|^q-B|u_{j-1}|^q)\|_2\\
&\le 2^{p-1}pA\|L'\left((|w_j|^{p-1}+|w_{j-1}|^{p-1})|w_j-w_{j-1}|\right)\|_2\\
&\quad+2^{q-1}qB\|L'\left((|u_j|^{q-1}+|u_{j-1}|^{q-1})|u_j-u_{j-1}|\right)\|_2\\
&\le 2^{p-1}pAC(T+R)(\|w_j\|_2^{p-1}+\|w_{j-1}\|_2^{p-1})\|w_j-w_{j-1}\|_2\\
&\quad+2^{q-1}qBC(T+R)^2(\|u_j\|_1^{q-1}+\|u_{j-1}\|_1^{q-1})\|u_j-u_{j-1}\|_1\\
&\le 2^ppAC(3M\e)^{p-1}(T+R)\|w_j-w_{j-1}\|_2\\
&\quad+2^qqBC(3M\e)^{q-1}(T+R)^2\|u_j-u_{j-1}\|_1.
\end{array}
\]
Here we employ H\"older's inequality to obtain
\[
\begin{array}{ll}
\|L(|w_j|^{p-1}|w_j-w_{j-1}|)\|_1
&=\|L\left(\left||w_j|^{(p-1)/p}|w_j-w_{j-1}|^{1/p}\right|^p\right)\|_1\\
&\le C(T+R)\||w_j|^{(p-1)/p}|w_j-w_{j-1}|^{1/p}\|_2^p\\
&\le C(T+R)\|w_j\|_2^{(p-1)}\|w_j-w_{j-1}\|_2
\end{array}
\]
and so on.
Therefore the convergence of $\{u_j\}$ follows from
\begin{equation}
\label{convergence}
\begin{array}{l}
\|u_{j+1}-u_j\|_1+\|w_{j+1}-w_j\|_2\\
\d\le\frac{1}{2}\left(\|u_j-u_{j-1}\|_1+\|w_j-w_{j-1}\|_2\right)
\end{array}
\quad(j\ge2)
\end{equation}
provided (\ref{condi1}) and
\begin{equation}
\label{condi2}
2^ppAC(3M\e)^{p-1}(T+R),2^qqBC(3M\e)^{q-1}(T+R)^2\le\frac{1}{4}
\end{equation}
are fulfilled.

\vskip10pt
\par\noindent
{\bf The convergence of the sequence $\v{\{\left((u_j)_x,(w_j)_x\right)\}}$.}
\par
First we note that $\|(u_1)_x\|_1,\|(w_1)_x\|_2\le M\e$ by (\ref{u^0_x}) and (\ref{u^0_tx}).
Assume that (\ref{condi1}) and (\ref{condi2}) are fulfilled.
Since (\ref{u_j_x,w_j_x}) and (\ref{apriori}) yield that
\[
\begin{array}{ll}
\|(u_{j+1})_x\|_1
&\le M\e+A\|\overline{L'}\left(|w_j|^p|\right)\|_1+B\|\overline{L'}\left(|u_j|^q|\right)\|_1\\
&\le M\e+AC(T+R)\|w_j\|_2^p+BC(T+R)^2\|u_j\|_1^q\\
&\le M\e+AC(3M\e)^p(T+R)+BC(3M\e)^q(T+R)^2
\end{array}
\]
because of a trivial property  $|\overline{L'}(v)|\le L'(|v|)$ and
\[
\begin{array}{ll}
\|(w_{j+1})_x\|_2
&\le M\e+pA\|L'\left(|w_j|^{p-1}|(w_j)_x|\right)\|_2\\
&\quad+qB\|L'\left(|u_j|^{q-1}|(u_j)_x|\right)\|_2\\
&\le M\e+pAC(T+R)\|w_j\|_2^{p-1}\|(w_j)_x\|_2\\
&\quad+qBC(T+R)^2\|u_j\|_1^{q-1}\|(u_j)_x\|_1\\
&\le M\e+pAC(3M\e)^{p-1}(T+R)\|(w_j)_x\|_2\\
&\quad+qBC(3M\e)^{q-1}(T+R)^2\|(u_j)_x\|_1.
\end{array}
\]
The boundedness of $\{\left((u_j)_x,(w_j)_x\right)\}$, i.e.
\begin{equation}
\label{bound_U_x}
\|(u_j)_x\|_1,\|(w_j)_x\|_2\le 3M\e\quad(j\in\N),
\end{equation}
follows from
\begin{equation}
\label{condi3}
pAC(3M\e)^p(T+R),qBC(3M\e)^q(T+R)^2\le M\e.
\end{equation}
Assuming (\ref{condi3}), one can estimate $\{(u_{j+1})_x-(u_j)_x\}$
and $\{(w_{j+1})_x-(w_j)_x\}$ as follows.
It is easy to see that
\[
\begin{array}{ll}
\|(u_{j+1})_x-(u_j)_x\|_1
&\le A\|\overline{L'}(|w_j|^p-|w_{j-1}|^p)\|_1\\
&\quad+B\|\overline{L'}(|u_j|^q-|u_{j-1}|^q)\|_1,
\end{array}
\]
which can be handled like $(w_{j+1}-w_j)$ as before, so that we have that
\[
\begin{array}{ll}
\|(u_{j+1})_x-(u_j)_x\|_1
&\le 2^ppAC(3M\e)^{p-1}(T+R)\|w_j-w_{j-1}\|_2\\
&\quad+2^qqBC(3M\e)^{q-1}(T+R)^2\|u_j-u_{j-1}\|_1
\end{array}
\]
because of $|\overline{L'}(v)|\le L'(|v|)$, which implies that
\begin{equation}
\label{convergence_u_x}
\|(u_{j+1})_x-(u_j)_x\|_1
=O\left(\frac{1}{2^j}\right)
\end{equation}
as $j\rightarrow\infty$ due to (\ref{convergence}).

\par
On the other hand, we have that
\[
\begin{array}{ll}
\|(w_{j+1})_x-(w_j)_x\|_2
&\le pA\|L'(|w_j|^{p-2}w_j(w_j)_x-|w_{j-1}|^{p-2}w_{j-1}(w_{j-1})_x)\|_2\\
&\quad+qB\|L'(|u_j|^{q-2}u_j(u_j)_x-|u_{j-1}|^{q-2}u_{j-1}(u_{j-1})_x)\|_2.
\end{array}
\]
The first term on the right hand side of this inequality is divided into two pieces
according to
\[
\begin{array}{l}
|w_j|^{p-2}w_j(w_j)_x-|w_{j-1}|^{p-2}w_{j-1}(w_{j-1})_x\\
=(|w_j|^{p-2}w_j-|w_{j-1}|^{p-2}w_{j-1})(w_j)_x\\
\quad+|w_{j-1}|^{p-2}w_{j-1}((w_j)_x-(w_{j-1})_x).
\end{array}
\]
Since one can employ the estimate
\[
\begin{array}{l}
\left||w_j|^{p-2}w_j-|w_{j-1}|^{p-2}w_{j-1}\right|
\\
\le
\left\{
\begin{array}{ll}
(p-1)2^{p-2}(|w_j|^{p-2}+|w_{j-1}|^{p-2})|w_j-w_{j-1}| & \mbox{for}\ p\ge2,\\
2|w_j-w_{j-1}|^{p-1} & \mbox{for}\ 1<p<2,
\end{array}
\right.
\end{array}
\]
and the same one in which $w$ is replaced with $u$, we obtain that
\[
\begin{array}{l}
\|(w_{j+1})_x-(w_j)_x\|_2\\
\le pAC(T+R)\|(w_j)_x\|_2\times\\
\qquad\times
\left\{
\begin{array}{ll}
(p-1)2^{p-2}(\|w_j\|_2^{p-2}+\|w_{j-1}\|_2^{p-2})\|w_j-w_{j-1}\|_2 & \mbox{for}\ p\ge2,\\
2\|w_j-w_{j-1}\|_2^{p-1} & \mbox{for}\ 1<p<2
\end{array}
\right.
\\
\quad+pAC(T+R)\|w_{j-1}\|_2^{p-1}\|(w_j)_x-(w_{j-1})_x\|_2\\
\quad +qBC(T+R)^2\|(u_j)_x\|_1\times\\
\qquad\times
\left\{
\begin{array}{ll}
(q-1)2^{q-2}(\|u_j\|_1^{q-2}+\|u_{j-1}\|_1^{q-2})\|u_j-u_{j-1}\|_1 & \mbox{for}\ q\ge2,\\
2\|u_j-u_{j-1}\|_1^{q-1} & \mbox{for}\ 1<q<2
\end{array}
\right.
\\
\quad+qBC(T+R)^2\|u_{j-1}\|_1^{q-1}\|(u_j)_x-(u_{j-1})_x\|_1.
\end{array}
\]
Hence it follows from (\ref{convergence}) and (\ref{convergence_u_x}) that
\[
\begin{array}{ll}
\|(w_{j+1})_x-(w_j)_x\|_2
&\le pAC(3M\e)^{p-1}(T+R)\|(w_j)_x-(w_{j-1})_x\|_2\\
&\quad\d+
O\left(\frac{1}{2^{j\min\{p-1,q-1,1\}}}\right) 
\end{array}
\]
as $j\rightarrow\infty$.
Therefore we obtain the convergence of $\{\left((u_j)_x,(w_j)_x\right)\}$ provided
\begin{equation}
\label{condi4}
pAC(3M\e)^{p-1}(T+R)\le\frac{1}{2}.
\end{equation}

\vskip10pt
\par\noindent
{\bf Continuation of the proof.}
\par
The convergence of the sequence $\{(u_j,w_j)\}$ to $(u,w)$ in the closed subspace of $X$
satisfying
\[
 \|u\|_1,\|(u_x)\|_1,\|w\|_2,\|(w)_x\|_2\le 3M\e
 \]
is established by (\ref{condi1}), (\ref{condi2}), (\ref{condi3}), and (\ref{condi4}),
which follow from
\[
C_0\e^{p-1}(T+R)\le1\quad\mbox{and}\quad C_0\e^{q-1}(T+R)^2\le1,
\]
where
\[
\begin{array}{ll}
C_0:=\max &
\{3^pACM^{p-1},3^qBCM^{q-1},\\
&\quad2^{p+2}3^{p-1}pACM^{p-1}, 2^{q+2}3^{q-1}qBCM^{q-1},\\
&\quad 3^ppACM^{p-1},3^3qBCM^{q-1},2\cdot3^{p-1}ACM^{p-1}\}.
\end{array}
\]
Therefore the statement of Theorem \ref{thm:lower-bound_non-zero}
is established with
\[
\left\{
\begin{array}{l}
\d c=\frac{1}{2}\min\{C_0^{-1},C_0^{-1/2}\},\\
\e_1:=\min\{(2C_0R)^{-1/(p-1)},(2^2C_0R^2)^{-2/(q-1)}\}
\end{array}
\right.
\]
because
\[
R\le\frac{1}{2}\min\{C_0^{-1}\e^{-(p-1)},C_0^{-1/2}\e^{-(q-1)/2}\}
\]
holds for $0<\e\le\e_1$.
\hfill$\Box$


\section{Proof of Proposition \ref{prop:apriori}}
\par
In this section, we prove a priori estimate (\ref{apriori}).
Recall the definition of $L$ in (\ref{nonlinear}) and $L'$ in (\ref{nonlinear_derivative}).
From now on, a positive constant $C$ independent of $T$ and $\e$
may change from line to line.

\par
It follows from the assumption on the supports and the definition of $L$ that
\[
|L(|w|^p)(x,t)|\le C\|w\|_2^pI_1(x,t)
\quad\mbox{for}\ |x|\le t+R,
\]
where we set
\[
I_1(x,t):=\int_0^tds\int_{x-t+s}^{x+t-s}(s-|y|+2R)^{-p}\chi_w(y,s)dy.
\]
Here, $\chi_w$ denotes a characteristic function of $\mbox{supp}\ w$.
First, we consider the case of $x\ge0$.
From now on, we employ the change of variables by
\[
\alpha=s+y,\beta=s-y.
\]
For $t+x\ge R$ and $t-x\ge R$, making use of  the symmetry of the weight in $y$ and
extending the domain of the integral, we have that
\[
\begin{array}{ll}
I_1(x,t)
&\d\le C\int_{-R}^Rd\beta\int_{-R}^{t+x}d\alpha\\
&\d\quad+C\int_R^{t-x}d\beta\int_{-R}^R(\alpha+2R)^{-p}d\alpha\\
&\d\quad+C\int_R^{t-x}(\beta+2R)^{-p}d\beta\int_R^{t+x}d\alpha\\
&\le C(t+x+R)\le C(T+R).
\end{array}
\]
For $t+x\ge R$ and $|t-x|\le R$, we also have that
\[
I_1(x,t)\le C\int_{-R}^{t-x}d\beta\int_{-R}^{t+x}d\alpha\le C(T+R).
\]
For $t+x\le R$, it is trivial that
\[
I_1(x,t)\le C.
\]
Summing up, we obtain that
\[
|L(|w|^p)(x,t)|\le C\|w\|_2^p(T+R)
\quad\mbox{for}\ 0\le x\le t+R.
\]
The case of $x\le0$ is similar to the one above, so we omit the details.
Therefore we obtain the first inequality of the first line of (\ref{apriori}).

\par
The second inequality of the first line of (\ref{apriori}) follows from
\[
|L(|u|^q)(x,t)|\le C\|u\|_1^qI_2(x,t)
\quad\mbox{for}\ |x|\le t+R,
\]
where
\[
I_2(x,t):=\int_0^tds\int_{x-t+s}^{x+t-s}\chi_u(y,s)dy.
\]
Here, $\chi_u$ denotes a characteristic function of $\mbox{supp}\ u$.
Indeed, it is trivial that
\[
I_2(x,t)\le C(T+R)^2\quad\mbox{for}\ |x|\le t+R.
\]

\par
Next, we shall show the second line in (\ref{apriori}).
It follows from the assumption on the supports and the definition of $L'$ that
\[
|L'(|w|^p)(x,t)|\le
C\|w\|_2^p\{I_+(x,t)+I_-(x,t)\}
\quad\mbox{for}\ |x|\le t+R,
\]
where the integrals $I_+$ and $I_-$ are defined by
\[
I_\pm(x,t):=\int_0^t(s-|t-s\pm x|+2R)^{-p}\chi_\pm(x,t;s)ds
\]
and the characteristic functions $\chi_+$ and $\chi_-$ are defined by
\[
\chi_\pm(x,t;s):=\chi_{\{s: |t-s\pm x|\le s+R\}},
\]
respectively.
First we note that it is sufficient to  estimate $I_\pm$ for $x\ge0$ due to its symmetry,
\[
I_+(-x,t)=I_-(x,t).
\]
Hence it follows from $0\le x\le t+R$ as well as
\[
|t-s+x|\le s+R\quad\mbox{and}\quad 0\le s\le t
\]
that
\[
\frac{t+x-R}{2}\le s\le t.
\]
Therefore we obtain
\[
I_+(x,t)\le\int_{(t+x-R)/2}^t(s-(t-s+x)+2R)^{-p}ds\le C
\quad\mbox{for}\ 0\le x\le t+R.
\]

\par
On the other hand, the estimate for $I_-$ is divided into two cases.
If $t-x\ge0$, then  $|t-s-x|\le s+R$ yields that
\[
\left\{
\begin{array}{ll}
t-s-x\le s+R & \mbox{for}\ 0\le s\le t-x,\\
-t+s+x\le s+R & \mbox{for}\ t-x\le s\le t,
\end{array}
\right.
\]
so that
\[
\begin{array}{ll}
I_-(x,t)
&\d\le\int_{(t-x-R)/2}^{t-x}(s-(t-s-x)+2R)^{-p}ds\\
&\d\quad+\int_{t-x}^t(s+(t-s-x)+2R)^{-p}ds\\
&\d\le C(1+(t-x+2R)^{-p}x)
\end{array}
\]
follows.
Therefore, neglecting $(t-x+2R)^{1-p}$, we obtain
\[
I_-(x,t)\le C(t-x+2R)^{-1}(T+R)\quad\mbox{for}\ t-x\ge0.
\]
If $(-R\le)t-x\le0$, $|t-s-x|\le s+R$ yields
\[
s-t+x\le s+R\quad\mbox{for}\quad 0\le s\le t,
\]
so that
\[
I_-(x,t)\le\int_0^t(s-(s-t+x)+2R)^{-p}ds.
\]
Therefore we obtain
\[
I_-(x,t)\le C(t-x+2R)^{-1}(T+R)\quad\mbox{for}\ -R\le t-x\le0.
\]
Summing up all the estimates for $I_+$ and $I_-$, we have
\[
|L'(|w|^p)(x,t)|\le
C\|w\|_2^p(t-|x|+2R)^{-1}(T+R)\quad\mbox{for}\ |x|\le t+R.
\]
This yields the first inequality of the second and third lines in (\ref{apriori}).

\par
The second inequality in the second and third lines in (\ref{apriori}) readily follows from
\[
|L'(|u|^q)(x,t)|\le
C\|u\|_1^q\{K_+(x,t)+K_-(x,t)\}
\quad\mbox{for}\ |x|\le t+R,
\]
where the integrals $K_+$ and $K_-$ are defined by
\[
K_\pm(x,t):=\int_0^t\chi_\pm(x,t;s)ds.
\]
Indeed, it is easy to obtain that
\[
K_\pm(x,t)\le C(t-|x|+2R)^{-1}(T+R)^2
\quad\mbox{for}\ |x|\le t+R.
\]
The proof of Proposition \ref{prop:apriori} is now completed.
\hfill$\Box$.


\section{Proof of Theorem \ref{thm:lower-bound_zero}}
\par
First we note that the strong Huygens' principle
\begin{equation}
\label{Huygens}
u^0(x,t)\equiv0\quad\mbox{in}\ D
\end{equation}
holds in this case of (\ref{zero}), where
\[
D:=\{(x,t)\in\R\times[0,\infty)\ :\ t-|x|\ge R\}.
\]
This is almost trivial if one takes a look on the representation of $u^0$ in (\ref{u^0})
and the support condition on the data in (\ref{supp_initial}).
But one can see also Proposition 2.2
in Kitamura, Morisawa and Takamura \cite{KMT22} for the details.
So, our unknown functions are
$U:=u-\e u^0$ and $W:=w-\e u_t^0$ in (\ref{system}).
Let $\{(U_j,W_j)\}_{j\in\N}$ be a sequence of $\{C^1(\R\times[0,T])\}^2$ defined by
\begin{equation}
\label{U_j,W_j}
\left\{
\begin{array}{ll}
U_{j+1}=L(A|W_j+\e u_t^0|^p+B|U_j+\e u^0|^q), & U_1=0,\\
W_{j+1}=L'(A|W_j+\e u_t^0|^p+B|U_j+\e u^0|^q), &  W_1=0.
\end{array}
\right.
\end{equation}
Then, $\left\{\left((U_j)_x,(W_j)_x\right)\right\}$ has to satisfy
\begin{equation}
\label{U_j_x,W_j_x}
\left\{
\begin{array}{ll}
(U_{j+1})_x&=\overline{L'}\left(A|W_j+\e u_t^0|^p+B|U_j+\e u^0|^q\right),\\
(U_1)_x &=0,\\
(W_{j+1})_x 
& =L'\left(Ap|W_j+\e u_t^0|^{p-2}(W_j+\e u_t^0)((W_j)_x+\e u_{tx}^0)\right)\\
& \quad+L'\left(Bq|U_j+\e u^0|^{q-2}(U_j+\e u^0)((U_j)_x+\e u_x^0)\right),\\
(W_1)_x &=0,
\end{array}
\right.
\end{equation}
so that the function space in which $\{(U_j,W_j)\}$ converges is
\[
\begin{array}{ll}
Y:=&\{(U,W)\in\{C^1(\R\times[0,T])\}^2\ :\ \|(U,W)\|_Y<\infty,\\
&\quad \mbox{supp}\ (U,W)\subset\{|x|\le t+R\}\}
\end{array}
\]
which is equipped with a norm 
\[
\|(U,W)\|_Y:=\|U\|_3+\|U_x\|_3+\|W\|_4+\|W_x\|_4,
\]
where
\[
\begin{array}{l}
\d\|U\|_3:=\sup_{(x,t)\in\R\times[0,T]}(t+|x|+R)^{-1}|U(x,t)|,\\
\d\|W\|_4:=\sup_{(x,t)\in\R\times[0,T]}\{\chi_D(x,t)+(1-\chi_D(x,t))(t+|x|+R)^{-1}\}|W(x,t)|,
\end{array}
\]
and $\chi_D$ is a characteristic function of $D$.
 Similarly to the proof of Theorem \ref{thm:lower-bound_non-zero},
we note that supp $(U_j,W_j)\in\{(x,t)\in\R\times[0,T]\ :\ |x|\le t+R\}$ implies
supp $(U_{j+1},W_{j+1})\in\{(x,t)\in\R\times[0,T]\ :\ |x|\le t+R\} $.

\par
The following lemmas are a priori estimates in this case.
\begin{prop}
\label{prop:apriori_linear}
Let $(U,W)\in\{C(\R\times[0,T])\}^2$ with
\[
\mbox{\rm supp}\ (U,W)\subset\{(x,t)\in\R\times[0,T]:|x|\le t+R\}
\]
and $U^0\in C(\R\times[0,T])$ with
\[
\mbox{\em supp}\ U^0\subset\{(x,t)\in\R\times[0,T]:(t-|x|)_+\le|x|\le t+R\}.
\]
Then there exists a positive constant $E$ independent of $T$ and $\e$ such that
\begin{equation}
\label{apriori_linear}
\left\{
\begin{array}{ll}
\|L(|U^0|^{p-m}|W|^m)\|_3 & \le E\|U^0\|_\infty^{p-m}\|W\|_4^m(T+R)^m,\\
\|L(|U^0|^{q-m}|U|^m)\|_3 & \le E\|U^0\|_\infty^{q-m}\|U\|_3^m(T+R)^m,\\
\|L'(|U^0|^{p-m}|W|^m)\|_4 & \le E\|U^0\|_\infty^{p-m}\|W\|_4^m(T+R)^m,\\
\|L'(|U^0|^{q-m}|U|^m)\|_4 & \le E\|U^0\|_\infty^{q-m}\|U\|_3^m(T+R)^m,\\
\|L'(|U^0|^{p-m}|W|^m)\|_3 & \le E\|U^0\|_\infty^{p-m}\|W\|_4^m(T+R)^m,\\
\|L'(|U^0|^{q-m}|U|^m)\|_3 & \le E\|U^0\|_\infty^{q-m}\|U\|_3^m(T+R)^m,
\end{array}
\right.
\end{equation}
where $p-m,q-m>0\ (m=0,1,2)$ and the norm $\|\cdot\|_\infty$ is defined by
\[
\|U^0\|_\infty:=\sup_{(x,t)\in\R\times[0,T]}|U^0(x,t)|.
\]
\end{prop}

\begin{prop}
\label{prop:apriori_zero}
Let $(U,W)\in\{C(\R\times[0,T])\}^2$ with
\[
\mbox{\rm supp}\ (U,W)\subset\{(x,t)\in\R\times[0,T]:|x|\le t+R\}.
\]
Then there exists a positive constant $C$ independent of $T$ and $\e$ such that
\begin{equation}
\label{apriori_zero}
\left\{
\begin{array}{ll}
\|L(|W|^p)\|_3 & \le C\|W\|_4^p(T+R)^p,\\
\|L(|U|^q)\|_3 & \le C\|U\|_3^q(T+R)^{q+1},\\
\|L'(|W|^p)\|_4 & \le C\|W\|_4^p(T+R)^p,\\
\|L'(|U|^q)\|_4 & \le C\|U\|_3^q(T+R)^{q+1},\\
\|L'(|W|^p)\|_3 & \le C\|W\|_4^p(T+R)^p,\\
\|L'(|U|^q)\|_3 & \le C\|U\|_3^q(T+R)^{q+1}.
\end{array}
\right.
\end{equation}
\end{prop}

The proofs of Proposition \ref{prop:apriori_linear} and \ref{prop:apriori_zero}
are established in the next section and after the next section respectively.
Set
\[
\begin{array}{ll}
N:=
&\d(2^ppA+2^qqB)E\left[\sum_{\gamma=0}^1\left\{\sum_{\alpha=0}^2
\left(\|f^{(\alpha)}\|_{L^\infty(\R)}^{p-\gamma}+\|f^{(\alpha)}\|_{L^\infty(\R)}^{q-\gamma}\right)\right.\right.\\
&\d\left.
+\|g\|_{L^1(\R)}^{p-\gamma}+\|g\|_{L^1(\R)}^{q-\gamma}+
\sum_{\beta=0}^1\left(\|g^{(\beta)}\|_{L^\infty(\R)}^{p-\gamma}
+\|g^{(\beta)}\|_{L^\infty(\R)}^{q-\gamma}\right)\right\}\\
&\d\left.+\sum_{\alpha=0}^2\|f^{(\alpha)}\|_{L^\infty(\R)}
+\|g\|_{L^1(\R)}+\sum_{\beta=0}^1\|g^{(\beta)}\|_{L^\infty(\R)}\right],
\end{array}
\]
where $E$ is the one in (\ref{apriori_linear}).
Since $u^0,u_t^0,u^0_x,u^0_{tx}$ satisfy
\[
\begin{array}{ll}
\|u^0\|_{\infty}
&\le\|f\|_{L^\infty(\mathbf{R})} + \|g\|_{L^1(\mathbf{R})},\\
\|u_t^0\|_{\infty},\ \|u^0_x\|_{\infty}
&\le \|f'\|_{L^\infty(\mathbf{R})}+\|g\|_{L^\infty(\mathbf{R})},\\
\|u^0_{tx}\|_{\infty}
&\le \|f''\|_{L^\infty(\mathbf{R})}+\|g'\|_{L^\infty(\mathbf{R})},
\end{array}
\]
we have that
\[
\begin{array}{ll}
\|u^0\|_{\infty}^p
&
\le 2^p\left(\|f\|_{L^\infty}^p+\|g\|_{L^1(\mathbf{R})}^p\right),\\
\|u^0_t\|_{\infty}^p,\ \|u^0_x\|_{\infty}^p
&\le 2^p\left(\|f'\|_{L^\infty(\mathbf{R})}^p+\|g\|_{L^\infty(\mathbf{R})}^p\right),\\
\|u^0_{tx}\|_{\infty}^p
&\le 2^p\left(\|f''\|_{L^\infty(\mathbf{R})}^p+\|g'\|_{L^\infty(\mathbf{R})}^p\right).
\end{array}
\]
Let us assume that
\[
0<\e\le1
\]
and define $\e_{2i}\ (i=1,2,3,4)$ respectively by
\begin{equation}
\label{epsilon}
\begin{array}{ll}
\e_{21}:= \min
&\left\{\left(2^{2p}3^pACN^{p-1}R^p\right)^{-1/(\min\{p,q\}(p-1))},\right.\\
&\quad\left.\left(2^{2q+1}3^qBCN^{q-1}R^{q+1}\right)^{-(q+1)/(\min\{p,q\}(q-1))}\right\},\\
\e_{22}:=\min
&\left\{\left(2^33^{2(p-1)}pACN^{p-1}R^p\right)^{-1/(\min\{p,q\}(p-1))},\right.\\
&\quad\left(2^23^{p-1}NR\right)^{-1/(p-1)},\\
&\quad\left(2^33^{2(q-1)}qBCN^{q-1}R^{q+1}\right)^{-(q+1)/(\min\{p,q\}(q-1))},\\
&\quad\left.\left(2^23^{q-1}NR\right)^{-1/(q-1)}\right\}\\
\e_{23}:=\min
&\left\{\e_{21},(2^{2p-1}3^{p+1}pACN^{p-1}R^p)^{-1/\min\{p,q\}(p-1)},\right.\\
&\quad\left(2^{p-1}3^pN^{p-1}R^{p-1} \right)^{-(p-1)/(\min\{p,q\}(p-2)+1)},\\
&\quad\left(2\cdot3^2NR \right)^{-1/(p-1)},\\
&\quad\left(2^{2q}3^{q+1}qBN^{q-1}R^{q+1}\right)^{-(q+1)/(\min\{p,q\}(q-1))},\\
&\quad\left(2^{q-1}3^{q}N^{q-1}R^{q-1} \right)^{-(q-1)/(\min\{p,q\}(q-2)+1)},\\
&\quad\left.\left(2\cdot3^2NR\right)^{-1/(q-1)}\right\},\\
\e_{24}:= \min
&\left\{\left(2^{2p+1}3^{p-1}pACN^{p-1}R^p\right)^{^{-1/(\min\{p,q\}(p-1))}},\right.\\
&\quad\left.\left(2^3NR\right)^{-1/(p-1)}\right\}.
\end{array}
\end{equation}

\vskip10pt
\par\noindent
{\bf The convergence of the sequence $\v{\{(U_j,W_j)\}}$.}
\par
It follows from (\ref{U_j,W_j}), Proposition \ref{prop:apriori_linear} and \ref{prop:apriori_zero}
that
\[
\begin{array}{ll}
\|U_{j+1}\|_3
&\le A\|L(|W_j + \e u^0_t|^p)\|_3 + B\|L(|U_j + \e u^0|^q)\|_3\\
&\le 2^pA\left\{\|L(|W_j|^p)\|_3+ \|L(\e |u^0_t|^p)\|_3\right\}\\
&\quad+2^q B\left\{\|L(|U_j|^q)\|_3+ \|L(\e|u^0|^q)\|_3 \right\}\\
&\le 2^pA\left\{C\|W_j\|_4^p(T+R)^p +E\e^p\|u_t^0\|_\infty^p\right\}\\
&\quad+2^q B\left\{C\|U_j\|_3^q(T+R)^{q+1} + E\e^q\|u^0\|_\infty^q \right\}.\\
&\le 2^pAC\|W_j\|_4^p(T+R)^{p}+2^qBC\|U_j\|_3^q(T+R)^{q+1} +N\e^{\min\{p,q\}}
\end{array}
\]
and
\[
\begin{array}{ll}
\|W_{j+1}\|_4
&\le A\|L'(|W_j + \e u_t^0 |^p)\|_4 + B\|L'(|U_j + \e u_0|^q)\|_4\\
&\le 2^pA\left\{\|L'(|W_j|^p)\|_4  + \|L'(\e| u^0_t|^p)\|_4 \right\}\\
&\quad+2^q B\left\{\|L'(|U_j|^q)\|_4 + \|L'(\e|u^0|^q)\|_4 \right\}\\
&\le 2^pA\left\{C\|W_j\|_4^p(T+R)^p + E\e^p\|u_t^0\|_\infty^p \right\}\\
&\quad+2^qB\left\{C\|U_j\|_3^q(T+R)^{q+1} + E\e^q\|u^0\|_\infty^q \right\}\\
&\le2^pAC\|W_j\|_4^p(T+R)^p +2^qBC\|U_j\|_3^{q}(T+R)^{q+1}+N\e^{\min\{p,q\}}.
\end{array}
\]
Hence the boundedness of $\{(U_j,W_j)\}$, i.e.
\begin{equation}
\label{bound_(U,W)}
\|U_j\|_3,\|W_j\|_4\le 3N\e^{\min\{p,q\}}\quad(j\in\N),
\end{equation}
follows from
\begin{equation}
\label{condi5}
\left\{
\begin{array}{ll}
2^pAC(3N\e)^{p\min\{p,q\}}(T+R)^p
&\le N\e^{\min\{p,q\}},\\
2^qBC(3N\e)^{q\min\{p,q\}}(T+R)^{q+1}
&\le N\e^{\min\{p,q\}}.
\end{array}
\right.
\end{equation}
Since (\ref{epsilon}) yields that
\[
R \le C_1\min\left\{\e^{-\min\{p,q\}(p-1)/p},\e^{-\min\{p,q\}(q-1)/(q+1)}\right\}
\]
for $0<\e\le\e_{21}$,
where
\[
C_1:=\frac{1}{2} \min\left\{\left(2^p3^{p}ACN^{p-1}\right)^{-1/p},
\left(2^q3^qBCN^{q-1}\right)^{-1/(q+1)}\right\},
\]
we find that (\ref{condi5}) as well as (\ref{bound_(U,W)}) follows from
\begin{equation}
\label{lifespan_step1}
T\le C_1\min\left\{\e^{-\min\{p,q\}(p-1)/p},\e^{-\min\{p,q\}(q-1)/(q+1)}\right\}
\end{equation}
for $0<\e\le\e_{21}$.

\par
Let us write down this inequality in each cases.
It follows from
\[
\frac{p-1}{p}-\frac{q-1}{q+1}= \frac{2p-q-1}{p(q+1)}
\]
that
\[
T\le C_1\e^{-p(p-1)/p} = C_1\e^{-(p-1)}
\quad\mbox{for}\ p\le \frac{q+1}{2}(<q)
\]
because of $\min\{p,q\} = p$,
\[
T\le C_1\e^{-p(q-1)/(q+1)}
\quad\mbox{for}\ \frac{q+1}{2}\le p\le q
\]
because of $\min\{p,q\} = p$, and
\[
T\le C_1\e^{-q(q-1)/(q+1)}
\quad\mbox{for}\ p\ge q
\]
because of $\min\{p,q\} = q$.

\par
Next, assuming (\ref{lifespan_step1}), one can estimate $(U_{j+1}-U_j)$ and $(W_{j+1}-W_j)$ as follows. The inequalities
\[
\begin{array}{l}
\|U_{j+1}-U_j\|_3\\
=\|L(A|W_j + \e u^0_t|^p-A|W_{j-1} + \e u^0_t|^p\\
\quad+
B|U_j + \e u^0|^q - B|U_{j-1} + \e u^0|^q)\|_3\\
\le  pA\|L(|W_{j-1} + \e u^0_t + \theta (W_j-W_{j-1})|^{p-1}|W_j-W_{j-1}|)\|_3\\
\quad+qB\|L(|U_{j-1} + \e u^0 + \theta(U_j-U_{j-1})|^{q-1}|U_j-U_{j-1}|)\|_3\\
\le
3^{p-1}pA\|L(|W_{j-1}|^{p-1} + |W_j|^{p-1} + \e^{p-1}|u^0_t|^{p-1})|W_j-W_{j-1}|\|_3\\
\quad+3^{q-1}qB\|L(|U_{j-1}|^{q-1} + |U_j|^{q-1} + \e^{q-1}|u^0|^{q-1})|U_j-U_{j-1}|\|_3 \\
\le 3^{p-1}pAC(\|W_{j-1}\|^{p-1}_4 + \|W_j\|^{p-1}_4)\|W_j-W_{j-1}\|_4(T+R)^p\\
\quad+3^{p-1}pAE\e^{p-1}\|u_t^0\|_\infty^{p-1}\|W_j-W_{j-1}\|_4(T+R)\\
\quad+3^{q-1}qBC(\|U_{j-1}\|^{q-1}_3 + \|U_j\|^{q-1}_3)\|U_j-U_{j-1}\|_3(T+R)^{q+1}\\
\quad+3^{q-1}qBE\e^{q-1}\|u^0\|^{q-1}_\infty\|U_j-U_{j-1}\|_3(T+R)\\
\le3^{p-1}\left\{2pAC(3N\e^{\min\{p,q\}})^{p-1}(T+R)^p
+N\e^{p-1}(T+R)\right\}\\
\qquad\times\|W_j-W_{j-1}\|_4\\
\quad+3^{q-1}\left\{2qBC(3N\e^{\min\{p,q\}})^{q-1}(T+R)^{q+1}
+N\e^{q-1}(T+R)\right\}\\
\qquad\times\|U_j-U_{j-1}\|_3
\end{array}
\]
and
\[
\begin{array}{l}
\|W_{j+1}-W_j\|_4\\
=\|L'(A|W_j + \e u_t^0 |^p -A|W_{j-1} + \e u_t^0 |^p\\
\quad+ B|U_j+ \e u_0|^q-B|U_{j-1} + \e u_0|^q)\|_4\\
\le 3^{p-1}pA\|L'(|W_{j-1}|^{p-1} + |W_j|^{p-1} + \e^{p-1}|u^0_t|^{p-1})|W_j-W_{j-1}|\|_4\\
\quad+3^{q-1}qB\|L'(|U_{j-1}|^{q-1} + |U_j|^{q-1} + \e^{q-1}|u^0|^{q-1})|U_j-U_{j-1}|\|_4\\
\le 3^{p-1}pAC(\|W_{j-1}\|^{p-1}_4 + \|W_j\|_4^{p-1})\|W_j-W_{j-1}\|_4(T+R)^p\\
\quad+3^{p-1}pAE\e^{p-1}\|u_t^0\|_\infty^{q-1}\|W_j-W_{j-1}\|_4(T+R)\\
\quad+3^{q-1}qBC(\|U_{j-1}\|_3^{q-1} + \|U_j\|_3^{q-1})\|U_j-U_{j-1}\|_3(T+R)^{q+1}\\
\quad+3^{q-1}qBE\e^{q-1}\|u^0\|_\infty^{q-1}\|U_j-U_{j-1}\|_3(T+R)\\
\le3^{p-1}\left\{2pAC(3N\e^{\min\{p,q\}})^{p-1})(T+R)^p
+N\e^{p-1}(T+R)\right\}\\
\qquad\times
\|W_j-W_{j-1}\|_4\\
\quad+3^{q-1}\left\{2qBC(3N\e^{\{\min\{p,q\}})^{q-1}(T+R)^{q+1}
+N\e^{q-1}(T+R))\right\}\\
\qquad\times\|U_j-U_{j-1}\|_3
\end{array}
\]
hold with some $\theta\in(0,1)$.
Here we employ H\"older's inequality like the one in the proof of Theorem
\ref{thm:lower-bound_non-zero}.
Therefore the convergence of $\{U_j\}$ follows from
\begin{equation}
\label{convergence_zero}
\begin{array}{l}
\|U_{j+1}-U_j\|_3+\|W_{j+1}-W_j\|_4\\
\d\le\frac{1}{2}\left(\|U_j-U_{j-1}\|_3+\|W_j-W_{j-1}\|_4\right)
\end{array}
\quad(j\ge2)
\end{equation}
provided (\ref{condi5}) and
\begin{equation}
\label{condi6}
\left\{
\begin{array}{ll}
2pAC(3N\e^{\min\{p,q\}})^{p-1}(T+R)^p+N\e^{p-1}(T+R)
&\d\le\frac{1}{3^{p-1}4},\\
2qBC(3N\e^{\min\{p,q\}})^{q-1}(T+R)^{q+1}+N\e^{q-1}(T+R)
&\d\le\frac{1}{3^{q-1}4}
\end{array}
\right.
\end{equation}
are fulfilled.
Since (\ref{epsilon}) yields that
\[
\begin{array}{ll}
R \le 
C_2\min
&\left\{\e^{-\min\{p,q\}(p-1)/p},\e^{-(p-1)}\right.,\\
&\quad\left.\e^{-\min\{p,q\}(q-1)/(q+1)},\e^{-(q-1)}\right\},
\end{array}
\]
for $0<\e\le\e_{22}$, where
\[
\begin{array}{ll}
C_2:=\d\frac{1}{2}\min
&\left\{
\left(2^33^{2(p-1)}pACN^{p-1}\right)^{-1/p},
(2^23^{p-1}pAN)^{-1},\right.\\
&\quad\left.\left(2^33^{2(q-1)}qBCN^{q-1}\right)^{-1/(q+1)},
(2^23^{q-1}qBN)^{-1}\right\},
\end{array}
\]
we find that (\ref{condi6}) as well as (\ref{convergence_zero}) follows from
\begin{equation}
\label{lifespan_step2}
T\le C_2
\min\left\{\e^{-\min\{p,q\}(p-1)/p},\e^{-(p-1)},
\e^{-\min\{p,q\}(q-1)/(q+1)},\e^{-(q-1)}\right\}
\end{equation}
for $0<\e\le\e_{22}$.

\par
It is trivial that
\[
\begin{array}{ll}
 \e^{-\min\{p,q\}(p-1)/p}=\e^{-(p-1)}\le\e^{-(q-1)} & \mbox{for}\ p\le q,\\
\e^{-\min\{p,q\}(q-1)/(q+1)\}}=\e^{-q(q-1)/(q+1)}\le\e^{-(q-1)} & \mbox{for}\ p\ge q.
\end{array}
\]
Therefore the computations after (\ref{lifespan_step1}) implies that
(\ref{lifespan_step2}) is equivalent to
\[
T\le
\left\{
\begin{array}{ll}
C_2\e^{-(p-1)} & \mbox{for}\ p\le\d\frac{q+1}{2}(<q),\\
C_2\e^{-p(q-1)/(q+1)} & \mbox{for}\ \d\frac{q+1}{2}\le p\le q,\\
C_2\e^{-q(q-1)/(q+1)} & \mbox{for}\ p\ge q.
\end{array}
\right.
\]

\vskip10pt
\par\noindent
{\bf The convergence of the sequence $\v{\{((U_j)_x,(W_j)_x)\}}$.}
\par
Assume (\ref{lifespan_step1}) and (\ref{lifespan_step2}).
Then we have (\ref{bound_(U,W)}) and (\ref{convergence_zero}).
Since it follows from (\ref{U_j_x,W_j_x}) that
\[
\begin{array}{ll}
|(U_{j+1})_x|
&\le 2^pA\{L'(|W_j|^p) + \e^pL'(|u^0_t|^p) \}\\
&\quad +2^qB\{L'(|U_j|^q + \e^qL'(|u^0|^q)\},\\
|(W_{j+1})_x|
&\le 2^{p-1}pA\{L'(|W_j|^{p-1}|(W_j)_x|)+ \e L'(|W_j|^{p-1}|u_{tx}^0|)\\
&\quad+\e^{p-1}L'(|u^0_t|^{p-1}|(W_j)_x|)+\e^pL'(|u^0_t|^{p-1}|u_{tx}|)\}\\
&\quad+2^{q-1}qB\{L'(|U_j|^{q-1}|(U_j)_x|) + \e L'(|U_j|^{q-1}|u^0|)\\
&\quad+\e^{q-1}L'(|u^0|^{q-1}|(U_j)_x| + \e^q L'(|u^0|^{q-1}|u^0_x|))\},
\end{array}
\]
Proposition \ref{prop:apriori_linear} and \ref{prop:apriori_zero} yield that
\[
\begin{array}{ll}
\|(U_{j+1})_x\|_3 
&\le2^pA\left\{L'(|W_j|^p)\|_3+\e^p\|L'(|u^0_t|^p)\|_3\right\}\\
&\quad+2^qB\left\{\|L'(|U_j|^q)\|_3+\e^q\|L'(|u^0|^q)\|_3\right\}\\
&\le2^pA\left\{C\|W_j\|^p_4(T+R)^p+\e^pE\|u_t^0\|_\infty^p\right\}\\
&\quad+2^qB\left\{C\|U_j\|^q_3(T+R)^{q+1}+\e^qE\|u^0\|_\infty^q\right\}\\
&\le N\e^{\min\{p,q\}}+2^pAC\|W_j\|_4^p(T+R)^p\\
&\quad+2^qBC\|U_j\|_3^q(T+R)^{q+1}\\
\end{array}
\]and
\[
\begin{array}{ll}
\|(W_{j+1})_x\|_4
&\le 2^{p-1}pA\left\{\|{L'}(|W_j|^{p-1}|(W_j)_x|)\|_4+\e \|L'(|W_j|^{p-1}|u^0_{tx}|)\|_3\right.\\
&\quad\left.+\e^{p-1}\|L'(|u^0_t|^{p-1}|(W_j)_x|)\|_4+\e^p\|L'(|u^0_t|^{p-1}|u_{tx}^0|)\|_4\right\}\\
&\quad+2^{q-1}qB\left\{\|L'(|U_j|^{q-1}|(U_j)_x|)\|_4+\e \|L'(|U_j|^{q-1}|u^0|)\|_4\right.\\
&\quad\left.+\e^{q-1}\|L'(|u^0|^{q-1}|(U_j)_x|)\|_4+\e^q\|L'(|u^0|^{q-1}|u^0_x|)\|_4\right\}\\
&\le 2^{p-1}pA\left\{C\|W_j\|_4^{p-1}\|(W_j)_x\|_4(T+R)^p\right.\\
&\qquad+\e E\|u^0_{tx}\|_\infty\|W_j\|_4^{p-1}(T+R)^{p-1}\\
&\quad\left.+\e^{p-1}E\|u^0_t\|_\infty^{p-1}\|(W_j)_x\|_4(T+R)
+\e^pE\|u^0_t\|_\infty^{p-1}\|u^0_{tx}\|_\infty\right\}\\
&\quad+2^{q-1}qB\left\{C\|U_j\|_3^{q-1}\|(U_j)_x\|_3(T+R)^{q+1}\right.\\
&\qquad\left.+\e E\|u^0\|_\infty\|U_j\|_3^{q-1}(T+R)^{q-1}\right.\\
&\quad\left.+\e^{q-1}E\|u^0\|_\infty^{q-1}\|(U_j)_x\|_3(T+R) 
+\e^qE\|u^0\|_\infty^{q-1}\|u^0_x\|_\infty\right\}\\
&\le N\e^{\min\{p,q\}}+2^{p-1}pAC(3N\e^{\min\{p,q\}})^{p-1}(T+R)^p\|(W_j)_x\|_4\\
&\quad+N\e (3N\e^{\min\{p,q\}})^{p-1}(T+R)^{p-1}+N\e^{p-1}(T+R)\|(W_j)_x\|_4\\
&\quad+2^{q-1}qBC(3N\e^{\min\{p,q\}})^{q-1}(T+R)^{q+1}\|(U_j)_x\|_3\\
&\quad+N\e(3N\e^{\min\{p,q\}})^{q-1}(T+R)^{q-1}+N\e^{q-1}(T+R)\|(U_j)_x\|_3.
\end{array}
\]
Hence the boundedness of $\left\{\left((U_j)_x,(W_j)_x\right)\right\}$, i.e.
\begin{equation}
\label{bound_(U_x,W_x)}
\|(U_j)_x\|_3,\|(W_j)_x\|_4
\le 
3N\e^{\min\{p,q\}}
\quad 
(j\in\mathbf{N}),
\end{equation}
follows from
\begin{equation}
\label{condi7}
\left\{
\begin{array}{ll}
2^pAC(3N\e^{\min\{p,q\}})^p(T+R)^p & \le N\e^{\min\{p,q\}},\\
2^qBC(3N\e^{\min\{p,q\}})^q(T+R)^{q+1} & \le N\e^{\min\{p,q\}}
\end{array}
\right.
\end{equation}
and
\begin{equation}
\label{condi8}
\left\{
\begin{array}{ll}
3\cdot2^{p-1}pAC(3N\e^{\min\{p,q\}})^{p}(T+R)^p & \le N\e^{\min\{p,q\}},\\
3N\e(3N\e^{\min\{p,q\}})^{p-1}(T+R)^{p-1} & \le N\e^{\min\{p,q\}},\\
3N\e^{p-1}(3N\e^{\min\{p,q\}})(T+R) & \le N\e^{\min\{p,q\}},\\
3\cdot2^{q-1}qB C(3N\e^{\min\{p,q\}})^{q}(T+R)^{q+1} & \le N\e^{\min\{p,q\}},\\
3N\e(3N\e^{\min\{p,q\}})^{q-1}(T+R)^{q-1} & \le N\e^{\min\{p,q\}},\\
3N\e^{q-1}(3N\e^{\min\{p,q\}})(T+R) & \le N\e^{\min\{p,q\}},\\
\end{array}
\right.
\end{equation}
Since (\ref{epsilon}) yields that
\[
\begin{array}{ll}
R\le C_3\min
&\left\{\e^{-\min\{p,q\}(p-1)/p},\e^{-(\min\{p,q\}(p-2)+1)/(p-1)},\e^{-(p-1)},\right.\\
&\quad\left.\e^{-\min\{p,q\}(q-1)/(q+1)},\e^{-(\min\{p,q\}(q-2) + 1)/(q-1)},\e^{-(q-1)}\right\}
\end{array}
\]
for $0<\e\le\e_{23}$, where
\[
\begin{array}{ll}
C_3:=\d\frac{1}{2}\min
&\left\{(2^p3^pACN^{p-1})^{-1/p},
(2^q3^qBCN^{q-1})^{-1/(q+1)},\right.\\
&\quad(2^{p-1}3^{p+1}pACN^{p-1})^{-1/p},
(3^pN^{p-1})^{-1/(p-1)},(3^2N)^{-1},\\
&\quad	
(2^{q-1}3^{q+1}qBCN^{q-1})^{-1/(q+1)},
\left.(3^qN^{q-1})^{-1/(q-1)},(3^2N)^{-1}\right\},
\end{array}
\]
we find that (\ref{condi7}) and (\ref{condi8}) as well as (\ref{bound_(U_x,W_x)}) follow from
\begin{equation}
\label{lifespan_step3}
\begin{array}{ll}
T\le C_3\min
&\left\{\e^{-\min\{p,q\}(p-1)/p},\e^{-(\min\{p,q\}(p-2)+1)/(p-1)},\right.\\
&\quad\e^{-(p-1)},\e^{-\min\{p,q\}(q-1)/(q+1)},\\
&\left.\quad\e^{-(\min\{p,q\}(q-2) + 1)/(q-1)},\e^{-(q-1)}\right\}
\end{array}
\end{equation}
for $0<\e\le\e_{23}$.

\par
Let us write down this inequality in each case.
For $p\le(q+1)/2(<q)$, we have that
\[
\begin{array}{l}
\d\frac{\min\{p,q\}(p-1)}{p} - (p-1)=0,\\
\d\frac{\min\{p,q\}(p-2) + 1}{p-1}-(p-1)=\frac{p(p-2) + 1}{p-1} - (p-1)=0,\\
\d\frac{\min\{p,q\}(q-1)}{q+1}-(p-1)=\frac{p(q-1)}{q+1} - (p-1)\\
\d\quad=\frac{-2p + q + 1}{q + 1}\ge\frac{-(q+1) + (q+1)}{q+1}=0,\\
\d\frac{\min\{p,q\}(q-2) + 1}{q-1} - (p-1)=\frac{p(q-2) + 1}{q-1} - (p-1)=\frac{q-p}{q-1} > 0,\\
(q-1) - (p-1) > 0
\end{array}
\]
which implies that (\ref{lifespan_step3}) is equivalent to
\[
T \le C_3 \e^{-(p-1)}\quad\mbox{for}\ p\le \frac{q + 1}{2}.
\]
For $(q + 1)/2\le p \le q$, we have that
\[
\begin{array}{l}
\d\frac{\min\{p,q\}(p-1)}{p} - \frac{p(q-1)}{q+1}=\frac{p(p-1)}{p} - \frac{p(q-1)}{q+1}\\
\d\quad= \frac{2p-q-1}{p(q+1)}\ge\frac{q+1-(q+1)}{p(q+1)}=0,\\
\d\frac{\min\{p,q\}(p-2) + 1}{p-1}-\frac{p(q-1)}{q+1}=\frac{p(p-2) + 1}{p-1} - \frac{p(q-1)}{q+1}\\
\d\quad=\frac{2p - q-1}{q + 1}\ge\frac{q + 1 - (q+1)}{q + 1}=0,\\
\d\frac{\min\{p,q\}(q-1)}{q+1} -\frac{p(q-1)}{q+1}=0,\\
\d\frac{\min\{p,q\}(q-2)+1}{q-1} - \frac{p(q-1)}{q+1}=\frac{p(q-2) + 1}{q-1} - \frac{p(q-1)}{q+1}\\
\d\quad=\frac{pq-3p+q+1}{(q-1)(q+1)}\ge\frac{p^2-3p+p+1}{(q-1)(q+1)}=\frac{(p-1)^2}{(q-1)(q+1)}>0,\\
\d(q-1)-\frac{p(q-1)}{q+1}\ge(q-1) - \frac{q(q-1)}{q+1}>0
\end{array}
\]
which implies that (\ref{lifespan_step3}) is equivalent to
\[
T\le C_3\e^{-p(q-1)/(q+1)}
\quad\mbox{for}\ \frac{q+1}{2}\le p \le q.
\]
For $p\ge q$, we have that
\[
\begin{array}{l}
\d\frac{\min\{p,q\}(p-1)}{p} - \frac{q(q-1)}{q+1}=\frac{q(p-1)}{p} - \frac{q(q-1)}{q+1}\\
\d\quad=\frac{2p-q-1}{p(q+1)}\ge\frac{q-1}{p(q+1)}>0,\\
\d\frac{\min\{p,q\}(p-2) + 1}{p-1}-\frac{q(q-1)}{q+1}=\frac{q(p-2) + 1}{p-1} - \frac{q(q-1)}{q+1}\\
\d\quad=\frac{2pq-q^2-2q+1}{(p-1)(q+1)}\ge\frac{(q-1)^2}{(p-1)(q+1)} >0,\\
\d\frac{\min\{p,q\}(q-1)}{q+1}-\frac{q(q-1)}{q+1}=0,\\
\d\frac{\min\{p,q\}(q-2) + 1}{q-1}-\frac{q(q-1)}{q+1}=\frac{q(q-2) + 1}{q-1} - \frac{q(q-1)}{q+1}\\
\d\quad=\frac{(q+1)(q-1)-q(q-1)}{q+1}> 0,\\
\d(q-1) -  \frac{q(q-1)}{q+1} > 0,
\end{array}
\]
which implies that (\ref{lifespan_step3}) is equivalent to
\[
T\le C\e^{-q(q-1)/(q+1)}
\quad\mbox{for}\ p\ge q.
\]

\par
Next, assuming (\ref{lifespan_step1}), (\ref{lifespan_step2}) and (\ref{lifespan_step3}),
we shall estimate $\{(U_{j+1})_x-(U_j)_x\}$ and $\{(W_{j+1})_x-(W_j)_x\}$.
It is easy to see that
\[
\begin{array}{l}
|(U_{j+1})_x - (U_j)_x|\\
\le L'(|A|W_j + \e u^0_t|^p - |W_{j-1} + \e u^0_t|^p|)\\
\quad+L'(|B|U_j + \e u^0|^{q} - |U_{j-1} + \e u^0|^q|)\\
\le pAL'(|W_{j-1} + \e u^0_t + \theta(W_{j+1} - W_j)|^{p-1}|W_j-W_{j-1}|)\\
\quad+qBL'(|U_{j-1} + \e u^0_t + \theta(U_{j+1}-U_j)|^{q-1}|U_j-U_{j-1}|)\\
\le 3^{p-1}pAL'\{(|W_{j-1}|^{p-1} + |W_j|^{p-1} + \e^{p-1}|u^0_t|^{p-1})|W_j-W_{j-1}|\}\\
\quad+3^{q-1}qBL'\{(|U_{j-1}|^{q-1} + |U_j|^{q-1} + \e^{q-1}|u^0|^{q-1})|U_j - U_{j-1}|\}
\end{array}
\]
hold for some $\theta\in(0,1)$.
Moreover we have that
\[
\begin{array}{l}
|(W_{j+1})_x - (W_j)_x|\\
\le L'(|pA|W_j + \e u^0_t|^{p-2}(W_j + \e u^0_t)((W_j)_x + \e u^0_{tx})\\
\quad-pA|W_{j-1} + \e u^0_t|^{p-2}(W_{j-1} + \e u^0_t)((W_{j-1})_x + \e u^0_{tx})|)\\
\quad+L'(qB|U_j + \e u^0|^{q-2}(U_j + \e u^0)((U_j)_x + \e u^0_x)\\
\quad-qB|U_{j-1} + \e u^0|^{q-2}(U_{j-1} + \e u^0)((U_{j-1})_x + \e u^0_x)|).
\end{array}
\]
In order to estimate the quantities in the right hand side of this inequality, we employ
\[
\begin{array}{l}
|W_j + \e u^0_t|^{p-2}(W_j + \e u^0_t)(W_j + \e u^0_t)_x\\
-|W_{j-1} + \e u^0_t|^{p-2}(W_{j-1} + \e u^0_t)(W_{j-1} + \e u^0_t)_x\\
=(|W_j + \e u^0_t|^{p-2}(W_j + \e u^0_t)
-|W_{j-1} + \e u^0_t|^{p-2}(W_{j-1} + \e u^0_t))(W_j + \e u^0_t)_x\\
\quad+|W_{j-1} + \e u^0_t|^{p-2}(W_{j-1} + \e u^0_t)((W_j + \e u^0_t)_x - (W_{j-1} + \e u^0_t)_x)
\end{array}
\]
and
\[
\begin{array}{l}
||W_j + \e u^0_t|^{p-2}(W_j + \e u^0_t)- |W_{j-1} + \e u^0_t|^{p-2}(W_{j-1} + \e u^0_t)|\\
\le\left\{
\begin{array}{ll}
(p-1)|W_j + \e u^0_t+ \theta (W_j-W_{j-1})|^{p-2}|W_j - W_{j-1}|
& \mbox{for}\ p\ge 2,\\
2|W_j - W_{j-1}|^{p-1} & \mbox{for}\ 1< p < 2,
\end{array}
\right.
\\
\le\left\{
\begin{array}{ll}
3^{p-2}(p-1)(|W_j|^{p-2}+|W_{j-1}|^{p-2} + |\e u^0_t|^{p-2})|W_j - W_{j-1}|
& \mbox{for}\ p\ge 2,\\
2|W_j - W_{j-1}|^{p-1} & \mbox{for}\ 1< p < 2,
\end{array}
\right.
\end{array}
\]
with some $\theta\in(0,1)$.
Hence it follows from Propositions \ref{prop:apriori_linear} and  \ref{prop:apriori_zero} that
\[
\begin{array}{l}
\|(U_{j+1})_x - (U_j)_x\|_3\\
\le 3^{p-1}pA\left\{C(\|W_{j-1}\|_4^{p-1} + \|W_j\|_4^{p-1})\|W_j - W_{j-1}\|_4(T+R)^p \right.\\
\quad+ \left.\e^{p-1}E\|u^0_t\|_\infty^{p-1}\|W_j - W_{j-1}\|_4(T+R)\right\}\\
\quad+3^{q-1}qB\left\{C(\|U_{j-1}\|_3^{q-1} + \|U_j\|_3^{q-1} )\|U_j - U_{j-1}\|_3(T+R)^{q+1}\right.\\
\quad+\left.\e^{q-1}E\|u^0\|_\infty^{q-1}\|U_j - U_{j-1}\|_3(T+R)\right\}
\end{array}
\]
and
\[
\begin{array}{l}
\|(W_j)_x - (W_{j-1})_x\|_4\\
\le\left\{
\begin{array}{ll}
\left.
\begin{array}{l}
3^{p-2}p(p-1)A\left\{C(\|W_{j-1}\|_4^{p-2} + \|W_j\|_4^{p-2})\right.\\
\times\|W_j -W_{j-1}\|_4\|(W_j)_x + \e (u^0_t)_x\|_4(T+R)^p\\
+\e^{p-2}E\|u_t^0\|_\infty^{p-2}\|W_j-W_{j-1}\|_4\\
\left.\times\|(W_j)_x + \e u^0_{tx}\|_4(T+R)^2\right\}
\end{array}
\right\}
& \mbox{for}\ p\ge 2,\\
\left.
\begin{array}{l}
2pAC\|W_j-W_{j-1}\|_4^{p-1}\\
\times\|(W_j)_x + \e u^0_{tx}\|_4(T+R)^p
\end{array}
\right\}
& \mbox{for}\ 1<p<2,
\end{array}\right.
\\
\quad+pA\|L'(|W_{j-1} + \e u^0_t|^{p-1}|(W_j)_x - (W_{j-1})_x|)\|_4\\
\quad+\left\{
\begin{array}{ll}
\left.
\begin{array}{l}
3^{q-2}q(q-1)B\left\{C(\|U_{j-1}\|_3^{q-2} + \|U_j\|_3^{q-2})\right.\\
\times \|U_j - U_{j-1}\|_3\|(U_j)_x+\e u^0_x\|_3(T+R)^{q+1}\\
+\e^{q-2}E\|u^0\|_\infty^{q-2}\|U_j-U_{j-1}\|_3\\
\left.\times\|(U_j)_x + \e u^0_x\|_3(T+R)^2\right\}
\end{array}
\right\}
& \mbox{for}\ q\ge 2,\\
\left.
\begin{array}{l}
2qBC\|U_j - U_{j-1}\|_3^{q-1}\\
\times\|(U_j)_x + \e u^0_x\|_3(T+R)^{q+1} 
\end{array}
\right\}
& \mbox{for}\ 1< q< 2,
\end{array}
\right.
\\
\quad+qBC\|U_{j-1}  + \e u^0\|_3^{q-1}\|(U_j)_x-(U_{j-1})_x\|_3(T+R)^{q+1}.
\end{array}
\]

\par
Since (\ref{lifespan_step1}) and (\ref{lifespan_step2}) yield (\ref{convergence_zero}),
we have that
\[
\|U_{j+1}-U_j\|_3 + \|W_{j+1} - W_{j}\|_3\le O\left(\frac{1}{2^j}\right)
\quad\mbox{as}\ j\rightarrow\infty.
\]
This fact implies that
\begin{equation}
\label{convergence_U_x}
\|(U_{j+1})_x-(U_j)_x\|_3 = O\left(\frac{1}{2^j}\right)
\quad\mbox{as}\ j\rightarrow\infty.
\end{equation}
Moreover it follows from (\ref{lifespan_step1}) as well as (\ref{bound_(U,W)}) that
\[
\begin{array}{ll}
pA\|L'(|W_{j-1} + \e u^0_t|^{p-1}|(W_j)_x - (W_{j-1})_x|)\|_4
\\
\le 2^{p-1}pA\|L'\{(|W_{j-1}|^{p-1} + \e^{p-1}|u^0_t|^{p-1})|(W_j)_x - (W_{j-1})_x|\}\|_4\\
\le 2^{p-1}pAC\|W_j\|_4^{p-1}\|(W_j)_x-(W_{j-1})_x\|_4(T+R)^p\\
\quad+2^{p-1}pAE\e^{p-1}\|u_t^0\|_\infty^{p-1}\|(W_j)_x-(W_{j-1})_x\|_4(T+R)\\
\le 2^{p-1}pAC(3N\e^{\min\{p,q\}})^{p-1}(T+R)^p\|(W_j)_x-(W_{j-1})_x\|_4\\
\quad+N\e^{p-1}(T+R)\|(W_j)_x-(W_{j-1})_x\|_4.
\end{array}
\]
Hence we have that
\[
\begin{array}{l}
\|(W_{j+1})_x-(W_j)_x\|_4\\
\le 2^{p-1}3^{p-1}pACN^{p-1}\e^{\min\{p,q\}(p-1)}(T+R)^p\|(W_j)_x-(W_{j-1})_x\|_4\\
\quad
+N\e^{p-1}(T+R)\|(W_j)_x-(W_{j-1})_x\|_4
+\d O\left(\frac{1}{2^{j\min\{p-1,q-1,1\}}}\right).
\end{array}
\]
as $j\rightarrow\infty$.
Therefore we obtain that
\begin{equation}
\label{convergence_W_x}
\|(W_{j+1})_x-(W_j)_x\|_4\le O\left(\frac{1}{2^{j}}\right)
\quad\mbox{as}\ j\rightarrow\infty
\end{equation}
provided
\begin{equation}
\label{condi9}
\left\{
\begin{array}{ll}
2^{p-1}3^{p-1}pACN^{p-1}\e^{\min\{p,q\}(p-1)}(T+R)^p & \le\d\frac{1}{4},\\
N\e^{p-1}(T+R) & \le\d\frac{1}{4}
\end{array}
\right.
\end{equation}
holds.
Since (\ref{epsilon}) yields that
\[
R\le C_4\min\{\e^{-\min\{p,q\}(p-1)/p},\e^{-(p-1)}\}
\]
for $0<\e\le \e_{24}$, where
\[
C_4:=\frac{1}{2} \min\left\{(2^{p+1}3^{p-1}pACN^{p-1})^{-1/p},(2^2N)^{-1}\right\},
\]
we find that (\ref{condi9}) as well as (\ref{convergence_U_x}) and (\ref{convergence_W_x})
follows from
\begin{equation}
\label{lifespan_step4}
T\le C_4\min\{\e^{-\min\{p,q\}(p-1)/p},\e^{-(p-1)}\}
\end{equation}
for $0<\e\le \e_{24}$.
We note that (\ref{lifespan_step4}) is equivalent to
\[
T\le\left\{
\begin{array}{ll}
C_4\e^{-(p-1)} & \mbox{for}\ p\le q,\\
C_4\e^{-q(p-1)/p} & \mbox{for}\ q\le p.
\end{array}
\right.
\]
Due to the computations after (\ref{lifespan_step3}),
we have that
\[
\e^{-p(q-1)/(q+1)}\le\e^{-(p-1)}\quad \mbox{for}\ \frac{q+1}{2}\le p\le q
\]
and
\[
\e^{-q(q-1)/(q+1)}\le\e^{-q(p-1)/p}\quad\mbox{for}\ q\le p.
\]

\vskip10pt
\par\noindent
{\bf Continuation of the proof.}
\par
The convergence of the sequence $\{(U_j,W_j)\}$ to $(U,W)$ in the closed subspace of $Y$
satisfying
\[
 \|U\|_3,\|(U_x)\|_3,\|W\|_4,\|(W)_x\|_4\le 3N\e^{\min\{p,q\}}
 \]
is established by  (\ref{epsilon}), (\ref{lifespan_step1}), (\ref{lifespan_step2}),
(\ref{lifespan_step3}) and (\ref{lifespan_step4}).
Therefore the statement of Theorem \ref{thm:lower-bound_zero} is established with
\[
\left\{
\begin{array}{l}
\d c=\min\left\{C_1,C_2,C_3,C_4\right\},\\
\e_2=\min\{1,\e_{21},\e_{22},\e_{23},\e_{24}\}.
\end{array}
\right.
\]
\hfill$\Box$


\section{Proof of Proposition \ref{prop:apriori_linear}}
\par
In this section, we prove a priori estimate (\ref{apriori_linear})
in the  similar manner to the one for (\ref{apriori}).   
A positive constant $C$ independent of $T$ and $\e$
may change from line to line.

\par
It follows from the assumption on the supports that
\[
\begin{array}{l}
|L(|U^0|^{p-m}|W|^m)(x,t)|\\
\le C\|U^0\|_\infty^{p-m}\|W\|_4^m(T+R)^mJ_0(x,t)
\end{array}
\quad\mbox{for}\ |x|\le t+R,
\]
where we set
\[
J_0(x,t):=\int_0^tds\int_{x-t+s}^{x+t-s}\chi_{\mbox{\footnotesize supp}\ U^0}(y,s)dy.
\]
First, we consider the case of $x\ge0$.
From now on, we employ the change of variables
\[
\alpha=s+y,\beta=s-y.
\]
For $t+x\ge R$ and $t-x\ge R$, extending the domain of the integral, we have that
\[
\begin{array}{ll}
J_0(x,t)
&\d\le C\int_{-R}^Rd\beta\int_{-R}^{t+x}d\alpha
+C\int_{-R}^{t-x}d\beta\int_{-R}^Rd\alpha\\
&\le C(t+x+R).
\end{array}
\]
For $t+x\ge R$ and $|t-x|\le R$, we also have that
\[
J_0(x,t)\le C\int_{-R}^Rd\beta\int_{-R}^{t+x}d\alpha\le C(t+x+R).
\]
For $t+x\le R$, it is trivial that $J_0(x,t)\le C$.
Summing up, we obtain that
\[
\begin{array}{l}
|L(|U^0|^{p-m}|W|^m)(x,t)|\\
\le C\|U^0\|_\infty^{p-m}\|W\|_3^m(T+R)^m(t+x+R)
\end{array}\quad
\mbox{for}\ 0\le x\le t+R.
\]
The case of $x\le0$ is similar to the one above, so we omit the details.
Therefore we obtain the first inequality in (\ref{apriori_linear}).
The second inequality in (\ref{apriori_linear}) follows from the computation above
because we have that
\[
\begin{array}{l}
|L(|U^0|^{q-m}|U|^m)(x,t)|\\
\le C\|U^0\|_\infty^{q-m}\|U\|_3^m(T+R)^mJ_0(x,t)
\end{array}
\quad\mbox{for}\ |x|\le t+R.
\]

\par
Next, we shall show the third inequality in (\ref{apriori_linear}).
It follows from the assumption on the support and the definition of $L'$ that
\[
\begin{array}{l}
|L'(|U^0|^{p-m}|W|^m)(x,t)|\\
\le C\|U^0\|_\infty^{p-m}\|W\|_4^m(T+R)^m\{J_{1+}(x,t)+J_{1-}(x,t)\}
\end{array}
\quad\mbox{for}\ |x|\le t+R,
\]
where the integrals $J_{1+}$ and $J_{1-}$ are defined by
\[
J_{1\pm}(x,t):=\int_0^t\chi_{1\pm}(x,t;s)ds
\]
and the characteristic functions $\chi_+$ and $\chi_-$ are defined by
\[
\chi_{1\pm}(x,t;s):=\chi_{\{s: |s-|t-s\pm x||\le R\}},
\]
respectively.
First we note that it is sufficient to  estimate $J_{1\pm}$ for $x\ge0$ due to its symmetry,
\[
J_{1+}(-x,t)=J_{1-}(x,t).
\]
For $(x,t)\in D\cap\{x\ge0\}$, we have that
\[
J_{1\pm}(x,t)=\int_{(t\pm x-R)/2}^{(t\pm x+R)/2}ds\le C.
\]
On the other hand, for $t+x\ge R,x\ge0$ and $|t-x|\le R$, we have that
\[
J_{1+}(x,t)=\int_{(t+x-R)/2}^tds\le C
\]
and
\[
J_{1-}(x,t)=\int_{(t-x-R)/2}^tds\le C(t+x+R).
\]
It is trivial that, for $t+x\le R, x\ge0$, we have $J_{1\pm}(x,t)\le C$.
Summing up all the estimates above, we obtain that
\[
\begin{array}{l}
|L'(|U^0|^{p-m}|W|^m)(x,t)|\\
\le C\|U^0\|_\infty^{p-m}\|W\|_4^m(T+R)^m\\
\quad\times\{\chi_D(x,t)+(1-\chi_D(x,t))(t+|x|+R)\}
\end{array}
\quad\mbox{for}\ |x|\le t+R.
\]
Therefore the third inequality in (\ref{apriori_linear}) is established.

\par
Due to the assumptions on the supports as well as the definitions of the norms,
it is clear that the proofs of the three remaining inequalities in (\ref{apriori_linear})
easily follow from the computations above.
The proof of Proposition \ref{prop:apriori_linear} is now completed.
\hfill$\Box$.


\section{Proof of Proposition \ref{prop:apriori_zero}}
\par
In this section, we prove a priori estimate (\ref{apriori_zero}).
Again, a positive constant $C$ independent of $T$ and $\e$
may change from line to line.

\par
It follows from the assumption on the supports and the definition of $L$ that
\[
|L(|W|^p)(x,t)|\le C\|W\|_4^pJ_1(x,t)
\quad\mbox{for}\ |x|\le t+R,
\]
where we set
\[
\begin{array}{ll}
J_1(x,t):=&\d\int_0^tds\int_{x-t+s}^{x+t-s}
\{\chi_D(s,y)\\
&\d+(1-\chi_D(s,y))(s+|y|+R)^p\}\chi_{\mbox{\footnotesize supp}\ W}(s,y)dy.
\end{array}
\]
First, we consider the case of $x\ge0$.
From now on, we employ the change of variables
\[
\alpha=s+y,\beta=s-y.
\]
For $(x,t)\in D$, extending the domain of the integral, we have that
\[
\begin{array}{ll}
J_1(x,t)
&\d\le C\int_{-R}^Rd\beta\int_{-R}^{x+t}(\alpha+R)^pd\alpha\\
&\d\quad+C\int_R^{t-x}(\beta+R)^pd\beta\int_{-R}^Rd\alpha\\
&\d\quad+C\int_R^{t-x}d\beta\int_R^{t+x}d\alpha\\
&\le C(T+R)^p(t+x+R).
\end{array}
\]
For $t+x\ge R$ and $|t-x|\le R$, we also have that
\[
J_1(x,t)\le C+C\int_{-R}^{t-x}d\beta\int_R^{t+x}(\alpha+R)^pd\alpha\le C(T+R)^p(t+x+R).
\]
For $t+x\le R$, it is trivial that
\[
J_1(x,t)\le C.
\]
Summing up, we obtain that
\[
|L(|W|^p)(x,t)|\le C\|W\|_4^p(T+R)^p(t+x+R)
\quad\mbox{for}\ 0\le x\le t+R.
\]
The case of $x\le0$ is similar to the one above, so we omit the details.
Therefore we obtain the first inequality of the first line of (\ref{apriori_zero}).

\par
The second inequality in (\ref{apriori_zero}) follows from
\[
|L(|U|^q)(x,t)|\le C\|U\|_3^qJ_2(x,t)
\quad\mbox{for}\ |x|\le t+R,
\]
where
\[
J_2(x,t):=\int_0^tds\int_{x-t+s}^{x+t-s}(s+|y|+R)^q
\chi_{\mbox{\footnotesize supp}\ U}(y,s)dy.
\]
It is trivial that
\[
J_2(x,t)\le C(T+R)^{q+1}(t+|x|+R)\quad\mbox{for}\ |x|\le t+R.
\]

\par
Next, we shall show the third inequality in (\ref{apriori_zero}).
It follows from the assumption on the supports and the definition of $L'$ that
\[
|L'(|W|^p)(x,t)|\le
C\|W\|_4^p\{J_+(x,t)+J_-(x,t)\}
\quad\mbox{for}\ |x|\le t+R,
\]
where the integrals $J_+$ and $J_-$ are defined by
\[
\begin{array}{ll}
J_\pm(x,t):=&
\d\int_0^t
\{\chi_{2\pm}(x,t;s)\\
&+(1-\chi_{2\pm}(x,t;s))(s+|t-s\pm x|+R)^p\}\chi_\pm(x,t;s)ds
\end{array}
\]
and the characteristic functions $\chi_+,\chi_-,\chi_{2+}$ and $\chi_{2-}$ are defined by
\[
\begin{array}{l}
\chi_\pm(x,t;s):=\chi_{\{s: |t-s\pm x|\le s+R\}},\\
\chi_{2\pm}(x,t;s):=\chi_{\{s: s-|t-s\pm x|\ge R\}}
\end{array}
\]
respectively.
First we note that it is sufficient to  estimate $J_\pm$ for $x\ge0$ due to its symmetry,
\[
J_+(-x,t)=J_-(x,t).
\]
For $(x,t)\in D\cap\{x\ge0\}$, we have
\[
\begin{array}{ll}
J_+(x,t)
&\d\le C\int_{(t+x-R)/2}^{(t+x+R)/2}(t+x+R)^pds
+C\int_{(t+x+R)/2}^tds\\
&\le C(t+x+R)^p+C(t+x+R)\\
&\le C(T+R)^p
\end{array}
\]
and
\[
\begin{array}{ll}
J_-(x,t)
&\d\le C\int_{(t-x-R)/2}^{(t-x+R)/2}(s+|t-s-x|+R)^pds
+C\int_{(t-x+R)/2}^tds\\
&\le C(t+x+R)^p+C(t+x+R)\\
&\le C(T+R)^p.
\end{array}
\]
For $t+x\ge R$ and $|t-x|\le R$, we have
\[
J_+(x,t)
\le C\int_{(t+x-R)/2}^t(t+x+R)^pds
\le C(t+x+R)(T+R)^p
\]
and
\[
\begin{array}{ll}
J_-(x,t)
&\d\le C\int_0^t(s+|t-s-x|+R)^pds\\
&\d\le C\int_0^{t-x}(t-x+R)^pds
+C\int_{t-x}^t(2s-t+x+R)^pds\\
&\le C(t-x+R)^{p+1}+C(t+x+R)^{p+1}\\
&\le C(t+x+R)(T+R)^p.
\end{array}
\]
It is trivial that $J_{\pm}(x,t)\le C$ for $t+x\le R$.
Therefore we obtain the third inequality in (\ref{apriori_zero}).

\par
The fourth, fifth and sixth  inequalities in (\ref{apriori_zero})
readily follow from the computations above.
The proof of Proposition \ref{prop:apriori_zero} is now completed.
\hfill$\Box$


\section{Proofs of Theorem \ref{thm:upper-bound_non-zero}
and Theorem \ref{thm:upper-bound_zero}}
\par
This section is devoted to the blow-up part of (\ref{lifespan_non-zero}) and (\ref{lifespan_zero}).
The proof can be established by a contradiction argument to the assumption that
$(u,w)$ is a continuous solution in the time interval $[0,T]$
of the associated system of integral equations (\ref{system}).

\vskip10pt
\par\noindent
{\bf Proof of Theorem \ref{thm:upper-bound_non-zero}.}
\par
Recall that our assumption on the initial data is (\ref{positive_non-zero}).
If one neglects $A|u_t|^p$ in (\ref{IVP_combined}), one immediately obtains
the contradiction if $T$ satisfies
\[
T\ge C\e^{-(q-1)/2}
\]
with suitably small $\e$, where $C$ is a positive constant independent of $\e$,
by Theorem 1.1 in Zhou \cite{Zhou92} or Theorem 5.1 in Takamura \cite{Takamura15}
for the classical solution of (\ref{IVP_combined}).
This estimate is also available for the continuous solution of (\ref{system})
if one employs the proof of Theorem 2.3 with $a=-1$
in Kitamura, Morisawa and Takamura \cite{KMT22},
which is applied to
\[
u\ge\e u^0+BL(|u|^q).
\]

\par
On the other hand, if one neglects $B|u|^q$ in (\ref{IVP_combined}),
one immediately obtains the contradiction if $T$ satisfies
\[
T\ge C\e^{-(p-1)}
\]
with suitably small $\e$, where $C$ is a positive constant independent of $\e$,
by Theorem 1.1 in Zhou \cite{Zhou01} for the classical solution of (\ref{IVP_combined}).
This estimate is also available for the continuous solution of (\ref{system})
if one employs the proof of Theorem 2.2 with $a=-1$
in Kitamura, Morisawa and Takamura \cite{KMT},
which is applied to
\[
w\ge\e u_t^0+AL'(|w|^p).
\]

\par
Therefore we obtain the contradiction if $T$ satisfies
\[
T\ge\min\{C\e^{-(p-1)},C\e^{-(q-1)/2}\}
\]
with suitably small $\e$, where $C$ is a positive constant independent of $\e$,
which asserts the statement of Theorem \ref{thm:upper-bound_non-zero}.
\hfill$\Box$

\vskip10pt
\par\noindent
{\bf Proof of Theorem \ref{thm:upper-bound_zero}.}
\par
Recall that our assumption on the initial data is (\ref{positive_zero}).
If one neglects $A|u_t|^p$ in (\ref{IVP_combined}), one immediately obtains the contradiction
if $T$ satisfies
\begin{equation}
\label{lifespan1}
T\ge C\e^{-q(q-1)/(q+1)}
\end{equation}
with suitably small $\e$, where $C$ is a positive constant independent of $\e$,
by Theorem 1.2 in Zhou \cite{Zhou92} or Theorem 5.1 in Takamura \cite{Takamura15}
for the classical solution of (\ref{IVP_combined}).
This estimate is also available for the continuous solution of (\ref{system})
if one employs the proof of Theorem 2.4 with $a=-1$
in Kitamura, Morisawa and Takamura \cite{KMT22},
which is applied to
\[
u\ge\e u^0+BL(|u|^q).
\]

\par
On the other hand, one can obtain the contradiction if $T$ satisfies
\begin{equation}
\label{lifespan2}
T\ge C\e^{-(p-1)}
\end{equation}
with suitably small $\e$, where $C$ is a positive constant independent of $\e$,
for the continuous solution of (\ref{system}).
To see this, one has to employ
\[
w\ge\e u_t^0+AL'(|w|^p)
\]
for which the definition of $L'$ implies
\[
w(x,t)\ge \e u_t^0(x,t)+A\int_0^t|w(x-t+s,s)|^pds.
\]
Setting $x-t=R/2$, we have $x+t=2t+R/2$ so that
\[
w\left(t+\frac{R}{2},t\right)>C_5\e+A\int_{R/4}^t\left|w\left(s+\frac{R}{2},s\right)\right|^pds
\quad\mbox{for}\ t\ge\frac{R}{4}
\]
due to (\ref{u^0_t}), where we set
\[
C_5:=-\frac{1}{4}f'\left(\frac{R}{2}\right)>0.
\]
Hence the comparison argument with a function $z$ satisfying
\[
z(t)=C_5\e+A\int_{R/4}^t|z(s)|^pds
\quad\mbox{for}\ t\ge\frac{R}{4}
\]
yields that
\[
w\left(t+\frac{R}{2},t\right)>z(t)\quad\mbox{for}\ t\ge\frac{R}{4},
\]
so that the desired contradiction can be obtained by the same arguments
as in the proof of Theorem 2.2 in Kitamura, Morisawa and Takamura \cite{KMT}. 
In fact, if we assume that there exists a point
\[
t_0:=\inf\left\{t>\frac{R}{4}\ :\ w\left(t+\frac{R}{2},t\right)=z(t)\right\},
\]
then we obtain a contradiction,
\[
\begin{array}{ll}
0
&\d=w\left(t_0+\frac{R}{2},t_0\right)-z(t_0)\\
&\d=A\int_{R/4}^{t_0}
\left\{\left|w\left(s+\frac{R}{2},s\right)\right|^p-|z(s)|^p\right\}ds>0,
\end{array}
\]
by $w(R/4+R/2,R/4)>z(R/4)$ and the continuity of $w,z$.
Taking into account of the fact that $z$ satisfies
\[
\left\{
\begin{array}{l}
z'=A|z|^p\quad\mbox{in}\ [R/4,\infty),\\
z(R/4)=C_5\e,
\end{array}
\right.
\]
we have a contradiction if $T$ satisfies
\[
T\ge C\e^{-(p-1)}\ge\{\mbox{The blow-up time of $z$}\}
\quad\mbox{for}\ 0<\e\le\e_4,
\]
where $C$ is a positive constant independent of $\e$ and
$\e_4$ is a constant depending on $A,C_5,R,p$. 
In fact, the representation of $z$ is 
\[
z(t)=\left\{(C_5\e)^{-(p-1)}-(p-1)A\left(t-\frac{R}{4}\right)\right\}^{-1/(p-1)}
\]
so that $z(t)\rightarrow\infty$ as
\[
t\rightarrow\frac{C_5^{-(p-1)}}{(p-1)A}\e^{-(p-1)}+\frac{R}{4}.
\]
If one sets
\[
\e_4=\left(\frac{4C_5^{-(p-1)}}{(p-1)AR}\right)^{1/(p-1)},
\]
then one has a contradiction, $w(T+R/2,T)=\infty$, provided
\[
T\ge\frac{2C_5^{-(p-1)}}{(p-1)A}\e^{-(p-1)}\ge\frac{C_5^{-(p-1)}}{(p-1)A}\e^{-(p-1)}+\frac{R}{4}
\]
for $0<\e\le\e_4$.
Therefore it follows from (\ref{lifespan1}) and (\ref{lifespan2}) that we obtain the contradiction
if $T$ satisfies
\[
T\ge\min\{C\e^{-(p-1)},C\e^{-q(q-1)/(q+1)}\}
\]
with suitably small $\e$, where $C$ is a positive constant independent of $\e$.
This estimate asserts Theorem \ref{thm:upper-bound_zero}
for $p\le(q+1)/2$, or $p\ge q$.
See (\ref{lifespan_zero}) and Remark \ref{rem:relation}.

\par
Hence we have to improve this estimate in the case of $(q+1)/2\le p\le q$.
To do this, we start with
\[
u=\e u^0+L(A|w|^p+B|u|^q).
\]
Set
\[
D_+:=\{(x,t)\ :\ t-x\ge R,\ x\ge0\}.
\]
Then, changing variables by $\alpha=s+y, \beta=s-y$, we have that
\[
L(|w|^p)(x,t)\ge\frac{1}{4}\int_{-R}^0d\beta\int_R^{t+x}|w(s,y)|^pd\alpha
\quad\mbox{for}\ (x,t)\in D_+.
\]
Hence it follows from
\[
w(y,s)\ge \e u_t^0(y,s)=-\frac{\e}{2}f'(y-s)
\quad\mbox{for}\ s+y\ge R,\ 0<y-s<R
\]
that
\begin{equation}
\label{iteration}
\left\{
\begin{array}{l}
\d u(x,t)\ge\frac{B}{4}\int_R^{t-x}d\beta\int_R^{t+x}|u(s,y)|^qd\alpha,\\
u(x,t)\ge AC_6\e^p(t+x-R),
\end{array}
\right.
\quad\mbox{for}\ (x,t)\in D_+,
\end{equation}
where we set
\[
C_6:=\frac{1}{2^{p+2}}\int_{-R}^0|f'(-\beta)|^pd\beta>0.
\]

\par
From now on, we employ a routine iteration procedure.
Assume an estimate
\begin{equation}
\label{estimate}
u(x,t)\ge M_n(t+x-R)^{a_n}(t-x-R)^{b_n}
\quad\mbox{for}\ (x,t)\in D_+
\end{equation}
holds, where $a_n,b_n\ge0$ and $M_n>0$.
The sequences $\{a_n\},\{b_n\}$ and $\{M_n\}$ are defined later.
Then it follows from the first line in (\ref{iteration}) that
\[
\begin{array}{ll}
u(x,t)
&\d\ge\frac{BM_n^q}{4}\int_R^{t-x}(\beta-R)^{qb_n}d\beta\int_R^{t+x}(\alpha-R)^{qa_n}d\alpha\\
&\d\ge\frac{BM_n^q}{4(qa_n+1)(qb_n+1)}(t+x-R)^{qa_n+1}(t-x-R)^{qb_n+1}
\end{array}
\]
for $(x,t)\in D_+$.
Hence (\ref{estimate}) holds for all $n\in\N$ provided
\[
\left\{
\begin{array}{ll}
a_{n+1}=qa_n+1,\ a_1=1,\\
b_{n+1}=qb_n+1,\ b_1=0
\end{array}
\right.
\]
and
\[
M_{n+1}\le\frac{BM_n^q}{4(qa_n+1)(qb_n+1)},\ M_1=AC_6\e^p.
\]
It is easy to see that
\[
a_n=\frac{q^n-1}{q-1},\ b_n=\frac{q^{n-1}-1}{q-1}\quad(n\in\N),
\]
which implies
\[
(qa_n+1)(qb_n+1)\le(qa_n+1)^2=a_{n+1}^2
\le\frac{q^{2(n+1)}}{(q-1)^2}.
\]
Therefore $M_n$ should be defined by
\[
M_{n+1}=BC_7q^{-2(n+1)}M_n^q,\ M_1=AC_6\e^p,
\]
where we set
\[
C_7:=\frac{(q-1)^2}{4}>0,
\]
so that (\ref{estimate}) implies that
\begin{equation}
\label{lower-bound}
u(x,t)\ge C_8\{(t+x-R)(t-x-R)\}^{-1/(q-1)}\exp\left\{Z(x,t)q^{n-1}\right\}
\end{equation}
for $(x,t)\in D_+$, where
\[
\begin{array}{ll}
Z(x,t)
&\d:=\frac{1}{q-1}\log\{(t+x-R)^q(t-x-R)\}\\
&\d\quad+\frac{1}{q-1}\log(BC_7)-2S_q\log q+\log(AC_6\e^p),\\
C_8
&\d:=\exp\left(-\frac{1}{q-1}\log (BC_7)\right)>0.
\end{array}
\]
Indeed, $M_n$ satisfies
\[
\log M_{n+1}=\log(BC_7)-2(n+1)\log q+q\log M_n,
\]
which implies
\[
\begin{array}{ll}
\log M_{n+1}
&=(1+q+\cdots+q^{n-1})\log(BC_7)\\
&\quad-2\{n+1+qn+\cdots+q^{n-1}(n+1-n+1)\}\log q\\
&\quad+q^n\log M_1\\
&\d=\frac{q^n-1}{q-1}\log(BC_7)-2q^{n-1}\log q\sum_{j=0}^n\frac{j+2}{q^j}+q^n\log M_1\\
&\d\ge-\frac{1}{q-1}\log(BC_7)+q^n\left\{\frac{1}{q-1}\log(BC_7)-2S_q\log q+\log M_1\right\},
\end{array}
\]
where we set
\[
S_q:=\sum_{j=0}^\infty\frac{j+2}{q^{j+1}}<\infty.
\]

\par
In view of (\ref{lower-bound}), if there exists a point $(x_0,t_0)\in D_+$ such that
$Z(x_0,t_0)>0$, we have a contradiction $u(x_0,t_0)=\infty$ to the fact that
$(u,w)$ is a continuous solution on the time interval $[0,T]$ with $T\ge t_0$
of (\ref{system}) by letting $n\rightarrow\infty$.
Let us set $t_0=2x_0$ and $t_0\ge4R$.
Then, since we have
\[
(t_0+x_0-R)^q(t_0-x_0-R)\ge\frac{1}{4}\left(\frac{5}{4}\right)^qt_0^{q+1},
\]
$Z(x_0,t_0)>0$ follows from
\[
t_0^{q+1}>4\left(\frac{4}{5}\right)^q\frac{q^{2(q-1)S_q}}{A^{(q-1)}BC_6^{(q-1)}C_7}\e^{-p(q-1)}.
\]
Therefore this inequality completes the proof of Theorem \ref{thm:upper-bound_zero}
for $(q+1)/2\le p\le q$ with $\e_4>0$ satisfying
\[
\left\{4\left(\frac{4}{5}\right)^q
\frac{q^{2(q-1)S_q}}{A^{(q-1)}BC_6^{(q-1)}C_7}
\right\}^{1/(q+1)}
\e_4^{-p(q-1)/(q+1)}=4R.
\]
\hfill$\Box$

\section*{Acknowledgement}
\par
The main part of this work was completed
when the first author was in the master course of Mathematical Institute, Tohoku University
and the third author had a second affiliation with
Research Alliance Center of Mathematical Sciences, Tohoku University.
The second author is partially supported
by the Grant-in-Aid for Young Scientists (No.18K13447), 
Japan Society for the Promotion of Science.
The third author is partially supported
by the Grant-in-Aid for Scientific Research
(A) (No.22H00097) and (B) (No.18H01132), 
Japan Society for the Promotion of Science.
All the authors appreciate Prof. Kunio Hidano (Mie Univ. in Japan)
for his pointing out our trivial oversight on the regularity of the solution
for low powers of the nonlinear terms,
and the referee for providing a number of corrections regarding the English language.


\bibliographystyle{plain}

\par
E-mail address: morisawa.katsuaki@musashi.ed.jp
\par
E-mail address: t-sasaki@musashino-u.ac.jp
\par
E-mail address: hiroyuki.takamura.a1@tohoku.ac.jp
\end{document}